\pdfoutput=1
\RequirePackage{lineno}
\documentclass[12pt,american,english,prb,showkeys]{revtex4}
\usepackage[T1]{fontenc}
\usepackage[latin9]{inputenc}
\usepackage{color}
\usepackage{babel}
\usepackage{units}
\usepackage{amsmath}
\usepackage{amssymb}
\usepackage{graphicx}
\usepackage[unicode=true,pdfusetitle,
 bookmarks=true,bookmarksnumbered=false,bookmarksopen=false,
 breaklinks=true,pdfborder={0 0 0},backref=false,colorlinks=true]
 {hyperref}
\usepackage{breakurl}
 \usepackage{enumerate}

\makeatletter
\@ifundefined{textcolor}{}
{%
 \definecolor{BLACK}{gray}{0}
 \definecolor{WHITE}{gray}{1}
 \definecolor{RED}{rgb}{1,0,0}
 \definecolor{GREEN}{rgb}{0,1,0}
 \definecolor{BLUE}{rgb}{0,0,1}
 \definecolor{CYAN}{cmyk}{1,0,0,0}
 \definecolor{MAGENTA}{cmyk}{0,1,0,0}
 \definecolor{YELLOW}{cmyk}{0,0,1,0}
}

\@ifundefined{definecolor}
 {\usepackage{color}}{}
\usepackage{algorithmic}

\@ifundefined{showcaptionsetup}{}{%
 \PassOptionsToPackage{caption=false}{subfig}}
\usepackage{subfig}

\makeatother

\linenumbers

\begin{document}

\title{Superiorization: An optimization heuristic for medical physics}

\author{Gabor T. Herman}
\thanks{Author to whom correspondence should be addressed}
\email{gabortherman@yahoo.com}
\homepage{http://www.dig.cs.gc.cuny.edu/~gabor/index.html}
\affiliation{Department of Computer Science, The Graduate Center, City University
of New York, New York, NY 10016, USA}

\author{Edgar Gardu{\~n}o}
\affiliation{Departamento de Ciencias de la Computaci{\'o}n, Instituto de Investigaciones
en Matem{\'a}ticas Aplicadas y en Sistemas, Universidad Nacional Aut{\'o}noma
de M{\'e}xico, Cd. Universitaria, C.P. 04510, Mexico City, Mexico}

\author{Ran Davidi}
\affiliation{Department of Radiation Oncology, Stanford University, Stanford,
CA 94305, USA}

\author{Yair Censor}
\affiliation{Department of Mathematics, University of Haifa, Mt. Carmel, 31905
Haifa, Israel}

\begin{abstract}
\linenumbers
\ \textbf{Purpose}: To describe and mathematically validate the superiorization
methodology, which is a recently-developed heuristic approach to optimization,
and to discuss its applicability to medical physics problem formulations
that specify the desired solution (of physically given or otherwise
obtained constraints) by an optimization criterion.

\textbf{Methods}: The superiorization methodology is presented as
a heuristic solver for a large class of constrained optimization problems.
The constraints come from the desire to produce a solution that is
constraints-compatible, in the sense of meeting requirements provided
by physically or otherwise obtained constraints. The underlying idea
is that many iterative algorithms for finding such a solution are
perturbation resilient in the sense that, even if certain kinds of
changes are made at the end of each iterative step, the algorithm
still produces a constraints-compatible solution. This property is
exploited by using permitted changes to steer the algorithm to a solution
that is not only constraints-compatible, but is also desirable according
to a specified optimization criterion. The approach is very general,
it is applicable to many iterative procedures and optimization criteria
used in medical physics.

\textbf{Results}: The main practical contribution is a procedure for
automatically producing from any given iterative algorithm its superiorized
version, which will supply solutions that are superior according to
a given optimization criterion. It is shown that if the original iterative
algorithm satisfies certain mathematical conditions, then the output
of its superiorized version is guaranteed to be as constraints-compatible
as the output of the original algorithm, but it is superior to the
latter according to the optimization criterion. This intuitive description
is made precise in the paper and the stated claims are rigorously
proved. Superiorization is illustrated on simulated computerized tomography
data of a head cross-section and, in spite of its generality, superiorization
is shown to be competitive to an optimization algorithm that is specifically
designed to minimize total variation.

\textbf{Conclusions}: The range of applicability of superiorization
to constrained optimization problems is very large. Its major utility
is in the automatic nature of producing a superiorization algorithm
from an algorithm aimed at only constraints-compatibility; while non-heuristic
(exact) approaches need to be redesigned for a new optimization criterion.
Thus superiorization provides a quick route to algorithms for the
practical solution of constrained optimization problems.
\end{abstract}


\maketitle

\section{\label{sec:Introduction}Introduction}

Optimization is a tool that is used in many areas of Medical Physics.
Prime examples are radiation therapy treatment planning and tomographic
reconstruction, but there are others such as image registration. Some
well-cited classical publications on the topic are \cite{DEAS97a,EZZE96a,GUST94a,GUST95a,LESS01a,MANZ03a,PUGA00a,SHEP02a,STUD97a,WU00a,YU96a,ZHAN04a}
and some recent articles are \cite{ABDO10a,BART11a,CHEN10a,FIEG11a,FRED11a,HOLD10A,HOLD11a,KIM11a,Men10a,PENF10a,SIDK11a,STAB11a,YANG10a,ZHAN11a}.

In a typical medical physics application, one uses \textit{constrained
optimization}, where the constraints come from the desire to produce
a solution that is \textit{constraints-compatible}, in the sense of
meeting the requirements provided by physically or otherwise obtained
constraints. In radiation therapy treatment planning, the requirements
are usually in the form of constraints prescribed by the treatment
planner on the doses to be delivered at specific locations in the
body. These doses in turn depend on information provided by an imaging
instrument, typically a Magnetic Resonance Imaging (MRI) or a Computerized
Tomography (CT) scanner. In tomography, the constraints come from
the detector readings of the instrument. In such applications, it
is typically the case that a large number of solutions would be considered
good enough from the point of view of being constraints-compatible;
to a large extent, but not entirely, due to the fact that there is
uncertainty as to the exact nature of the constraints (for example,
due to noise in the data collection). In such a case, an optimization
criterion is introduced that helps us to distinguish the ``better''
constraints-compatible solutions (for example, this criterion could
be the total dose to be delivered to the body, which may vary quite
a bit between radiation therapy treatment plans that are compatible
with the constraints on the doses delivered to individual locations).

The superiorization methodology (see, for example, \cite{PENF10a,CENS10a,DAVI09a,NIKA12a,HERM08a,BUTN07a,GARD11a})
is a recently-developed heuristic approach to optimization. The word
\textit{heuristic} is used here in the sense that the process is not
guaranteed to lead to an optimum according to the given criterion;
approaches aimed at processes that are guaranteed in that sense are
usually referred to as \textit{exact.} Heuristic approaches have been
found useful in practical applications of optimization, mainly because
they are often computationally much less expensive than their exact
counterparts, but nevertheless provide solutions that are appropriate
for the application at hand \cite{RARD01a,WERN00a,ZANA81a}.

The underlying idea of the superiorization approach is the following.
In many applications there exists a computationally-efficient iterative
algorithm that produces a constraints-compatible solution for the
given constraints. (An example of this for radiation therapy treatment
planning is reported in \cite{Herman2008}, its clinical use is discussed
in \cite{CHEN10a}.) Furthermore, often the algorithm is \textit{perturbation
resilient} in the sense that, even if certain kinds of changes are
made at the end of each iterative step, the algorithm still produces
a constraints-compatible solution \cite{CENS10a,DAVI09a,NIKA12a,BUTN07a}.
This property is exploited in the superiorization approach by using
such perturbations to steer the algorithm to a solution that is not
only constraints-compatible, but is also desirable according to a
specified optimization criterion. The approach is very general, it
is applicable to many iterative procedures and optimization criteria.

The current paper presents a major advance in the practice and theory
of superiorization. The previous publications\cite{PENF10a,CENS10a,DAVI09a,NIKA12a,HERM08a,BUTN07a,GARD11a}
used the intuitive idea to present some superiorization algorithms,
in this paper the reader will find a totally automatic procedure that
turns an iterative algorithm into its superiorized version. This version
will produce an output that is as constraints-compatible as the output
of the original algorithm, but it is superior to that according to
an optimization criterion. This claim is mathematically shown to be
true for a very large class of iterative algorithms and for optimization
criteria in general, typical restrictions (such as convexity) on the
optimization criterion are not essential for the material presented
below. In order to make precise and validate this broad claim, we
present here a new theoretical framework. The framework of \cite{CENS10a}
is a precursor of what we present here, but it is a restricted one,
since it assumes that the constraints can be all satisfied simultaneously,
which is often false in medical physics applications. There is no
such restriction in the presentation below.

The idea of designing algorithms that use interlacing steps of two
different kinds (in our case, one kind of steps aim at constraints-compatibility
and the other kind of steps aim at improvement of the optimization
criterion) is well-established and, in fact, is made use of in many
approaches that have been proposed with exact constrained optimization
in mind; see, for example, the works of Helou Neto and De Pierro\cite{hdp-siam09,hdp11},
of Nurminski\cite{nur10}, of Combettes and coworkers\cite{COMB02a,COMB04a},
of Sidky and Pan and coworkers\cite{SIDK11a,SIDK08a,BIAN10A} and
of Defrise and coworkers\cite{DEFR11a}. However, none of these approaches
can do what can be done by the superiorization approach as presented
below, namely the automatic production of a heuristic constrained
optimization algorithm from an iterative algorithm for constraints-compatibility.
For example, in \cite{hdp-siam09} it is assumed (just as in the theory
presented in our \cite{CENS10a}) that all the constraints can be
satisfied simultaneously. 

A major motivator for the additional theory presented in the current
paper is to get rid of this assumption, which is not reasonable when
handling real problems of medical physics. Motivated by similar considerations,
Helou Neto and De Pierro \cite{hdp11} present an alternative approach
that does not require this unreasonable assumption. However, in order
to solve such a problem, they end up with iterative algorithms of
a particular form rather than having the generality of being able
to turn any constraints-compatibility seeking algorithm into a superiorized
one capable of handling constrained optimization. Also, the assumptions
they have to make in order to prove their convergence result (their
Theorem 15) indicate that their approach is applicable to a smaller
class of constrained optimization problems than the superiorization
approach whose applicability seems to be more general. However, for
the mathematical purist, we point out that they present an exact constrained
optimization algorithm, while superiorization is a heuristic approach.
Whether this is relevant to medical physics practice is not clear:
exact algorithms are not run forever, but are stopped according to
some stopping-rule, the relevant questions in comparing two algorithms
are the quality of the actual output and the computation time needed
to obtain it.

Ultimately, the quality of the outputs should be evaluated by some
figures of merit relevant to the medical task at hand. An example
of a careful study of this kind that involves superiorization is in
\cite[Section 4.3]{NIKA12a}, which reports on comparing in CT the
efficacy of constrained optimization reconstruction algorithms for
the detection of low-contrast brain tumors by using the method of
statistical hypothesis testing (which provides a P-value that indicates
the significance by which we can reject the null hypothesis that the
two algorithms are equally efficacious in favor of the alternative
that one is preferable). Such studies bundle together two things:
(i) the formulation of the constrained optimization task and (ii)
the performance of the algorithm in performing that task. The first
of these requires a translation of the medical aim into a mathematical
model, it is important that this model should be appropriately chosen. 

The superiorization approach is not about choosing this model, it
kicks in once the model is chosen and aims at producing an output
that is ``good'' according to the mathematical specifications of
the constraints and of the optimization criterion. Thus superiorization
has been used to compare the effects on the quality of the output
in CT when the optimization criterion is specified by total variation
(TV) versus by entropy\cite{DAVI09a} or versus by the $\ell_{1}$-norm
of the Haar transform \cite{GARD11a}. However, the current paper
is not about discussing how to translate the underlying medical physics
task into a constrained optimization problem. For our purposes here,
we are assuming that the mathematical model has been worked out and
concentrate on the algorithmic approach for solving the resulting
constrained optimization problem. We claim that the evaluation of
such algorithms should not be based on the medical figures of merit
mentioned at the beginning of the previous paragraph, but rather on
their performance in solving the mathematical problem. If ``good''
solutions to the constrained optimization problem are not medically
efficacious, that indicates that something is wrong with the mathematical
model and not that something is wrong with the algorithmic approach.
For this reason, in this paper we will not carry out a careful investigation
of the medical efficacy of any algorithm in the manner that we have
done in \cite[Section 4.3]{NIKA12a}, but will restrict ourselves
to a simple illustration of the performance of the superiorization
approach as compared to the previously published algorithm of\cite{SIDK08a}
that is aimed at performing exact minimization.

Examples of such studies already exist. Superiorization was compared
in \cite{BUTN07a} with Algorithm 6 of \cite{COMB02a} and in \cite{CENS12a}
with the algorithm of Goldstein and Osher that they refer to as TwIST
\cite{BIOU07a} with split Bregman \cite{GOLD09a} as the substep.
In both cases the implementation was done by the proposers of the
algorithms. In these reported instances superiorization did well:
the constraints-compatibility and the value of the function to be
minimized were very similar for the outputs produced by the algorithms
being compared, but the superiorization algorithm produced its output
four times faster than the alternative. It would be unjustified to
draw any general conclusions on the mathematical performance and speed
of superiorization based on just a few experiments, but the reported
results are encouraging.

However, the main reason why we advocate superiorization is different
from what is discussed above. The reason why we claim it to be helpful
in medical physics research is that it has the potential of saving
a lot of time and effort for the researcher. Let us consider a historical
example. Likelihood optimization using the iterative process of expectation
maximization (EM) \cite{SHEP82a} gained immediate and wide acceptance
in the emission tomography community. It was observed that irregular
high amplitude patterns occurred in the image with a large number
of iterations, but it was not until five years later that this problem
was corrected \cite{LEVI87A} by the use of a maximum a posteriority
probability (MAP) algorithm with a multivariate Gaussian prior. Had
we had at our disposal the superiorization approach, then the introduction
of an optimization criterion (Gaussian or other) into the iterative
expectation maximization (EM) process would have been a simple matter
and we would have saved the time and effort spent on designing a special
purpose algorithm for the MAP formulation. A $TV$-superiorization
of the EM algorithm is presented in\cite{JIN12a}.

Even though our major claim for superiorization is that it provides
a quick route to algorithms for the practical solution of constrained
optimization problems, before leaving this introduction let us bring
up a question that has to do with the performance of the resulting
algorithms: Will superiorization produce superior results to those
produced by contemporary MAP methods or is it faster than the
better of such methods? At this stage we have not yet developed
the mathematical notation to discuss this question in a rigorous manner,
we return to it in Subsection \ref{sub:Information-on-performance}.

In the next section we present in detail the superiorization methodology.
In the subsequent section we provide an illustrative example by reporting
on reconstructions produced by algorithms applied to simulated computerized
tomography data of a head cross-section. In the final section we discuss
our results and present our conclusions.

\section{The Superiorization Methodology\label{sec:Methods}}

\subsection{\label{sub:The-superiorization-methodology}Problem sets, proximity
functions and $\varepsilon$-compatibility\label{sub:superiorization}}

Although optimization is often studied in a more general context (such
as in Hilbert or Banach spaces), in medical physics we usually deal
with a special case, where optimization is performed in a \textit{Euclidean
space} $\mathbb{R}^{J}$ (the space of
$J$-dimensional vectors of real numbers, where $J$ is a positive
integer). As often appropriate in practice, we further restrict the
domain of optimization to a nonempty subset $\Omega$ of $\mathbb{R}^{J}$
(such as the \textit{nonnegative orthant} $\mathbb{R}_{+}^{J}$ that
consists of vectors all of whose components are nonnegative).

We now turn to formalizing the notion of being compatible with given
constraints, a notion that we have used informally in the previous
section. In any application, there is a \emph{problem set} $\mathbb{T}$;
each \textit{problem} $T\in\mathbb{T}$ is essentially a description
of the constraints in that particular case. For example, for a tomographic
scanner, the problem of reconstruction for a particular patient at
a particular time is determined by the measurements taken by the scanner
for that patient at that time. The intuitive notion of constraints-compatibility
is formalized by the use of a\textit{ proximity function} $\mathcal{P}r$
on $\mathbb{T}$ such that, for every $T\in\mathbb{T}$, $\mathcal{P}r_{T}$
maps $\Omega$ into $\mathbb{R}_{+}$, the set of nonnegative real
numbers; i.e., $\mathcal{P}r_{T}:\Omega\rightarrow\mathbb{R}_{+}$.
Intuitively we think of $\mathcal{P}r_{T}\left(\boldsymbol{x}\right)$
as an indicator of how incompatible $\boldsymbol{x}$ is with the
constraints of $T$. For example, in tomography, $\mathcal{P}r_{T}\left(\boldsymbol{x}\right)$
should indicate by how much a proposed reconstruction that is described
by an $\boldsymbol{x}$ in $\Omega$ violates the constraints of the
problem $T$ that are provided by the measurements taken by the scanner.
For example, if we use $\boldsymbol{b}$ to denote the vector of estimated
line integrals based on the measurements obtained by the scanner and
by $\boldsymbol{A}$ the system matrix of the scanner, then a possible
choice for the proximity function is the norm-distance $\left\Vert \boldsymbol{b-Ax}\right\Vert $,
which we will use as an example in the discussions that follow. An
alternative legitimate choice for the proximity function is the Kullback-Leibler
distance $KL(\boldsymbol{b},\boldsymbol{Ax})$, which is the negative
log-likelihood of a statistical model in tomography. The special
case $\mathcal{P}r_{T}\left(\boldsymbol{x}\right)=0$ is interpreted
by saying that $\boldsymbol{x}$ is perfectly compatible with the
constraints; due to the presence of noise in practical applications,
it is quite conceivable that there is no $\boldsymbol{x}$ that is
perfectly compatible with the constraints, and we accept an $\boldsymbol{x}$
as constraints-compatible as long as the value of $\mathcal{P}r_{T}\left(\boldsymbol{x}\right)$
is considered to be small enough to justify that decision. Combining
these two concepts leads to the notion of a \emph{problem} \emph{structure,}
which is a pair $\left\langle \mathbb{T},\mathcal{P}r\right\rangle $,
where $\mathbb{T}$ is a nonempty problem set and $\mathcal{P}r$
is a proximity function on $\mathbb{T}$.
For a problem structure $\left\langle \mathbb{T},\mathcal{P}r\right\rangle $,
a problem $T\in\mathbb{T}$, a nonnegative $\varepsilon$ and an $\boldsymbol{x}\in\Omega$,
we say that $\boldsymbol{x}$ is $\varepsilon$\textit{-compatible}
with $T$ provided that $\mathcal{P}r_{T}\left(\boldsymbol{x}\right)\leq\varepsilon$.

As an example (whose applicability to tomographic reconstruction is
illustrated in Section \ref{sec:An-Illustrative-ExampleResults}),
consider the problem structure that arises from the desire to find
nonnegative solutions of sequences of blocks of linear equations.
Then the appropriate choices are $\Omega=\mathbb{R}_{+}^{J}$ and
the problem structure is $\left\langle \mathbb{S},Res\right\rangle $,
where the problem set $\mathbb{S}$ is 
\begin{linenomath}
\begin{equation}
\begin{array}{lll}
\mathbb{S} & = & \left\{ \left(\left\{ \left(\boldsymbol{a}^{1},b_{1}\right),\ldots,\left(\boldsymbol{a}^{\ell_{1}},b_{\ell_{1}}\right)\right\} ,\ldots,\right.\right.\\
 &  & \left.\left.\left\{ \left(\boldsymbol{a}^{\ell_{1}+\ldots+\ell_{W-1}+1},b_{\ell_{1}+\ldots+\ell_{W-1}+1}\right),\ldots,\left(\boldsymbol{a}^{\ell_{1}+\ldots+\ell_{W}},b_{\ell_{1}+\ldots+\ell_{W}}\right)\right\} \right)\right|\\
 &  & W\mbox{ is a positive integer and,}\\
 &  & \mbox{for \ensuremath{1\leq w\leq W,\;\ell_{w}\:}is a positive integer and,}\\
 &  & \left.\mbox{for }\ensuremath{1\leq i\leq\ell_{1}+\ldots+\ell_{W},\:}\ensuremath{\boldsymbol{a}^{i}\in\mathbb{R}^{J}}\mbox{ and }b_{i}\in\mathbb{R}\right\} 
\end{array}\label{eq:S_definition}
\end{equation}
\end{linenomath}
and the proximity function $Res$ on $\mathbb{S}$ is defined, for
any problem $S=\left(\left\{ \left(\boldsymbol{a}^{1},b_{1}\right),\right.\right.$
$\left.\ldots,\left(\boldsymbol{a}^{\ell_{1}},b_{\ell_{1}}\right)\right\} ,\left.\ldots,\left\{ \left(\boldsymbol{a}^{\ell_{1}+\ldots+\ell_{W-1}+1},b_{\ell_{1}+\ldots+\ell_{W-1}+1}\right),\ldots,\left(\boldsymbol{a}^{\ell_{1}+\ldots+\ell_{W}},b_{\ell_{1}+\ldots+\ell_{W}}\right)\right\} \right)$
in $\mathbb{S}$ and for any $\boldsymbol{x}\in\Omega$, by 
\begin{linenomath}
\begin{equation}
Res_{S}(\boldsymbol{x})=\sqrt{\sum\limits _{i=1}^{\ell_{1}+\ldots+\ell_{W}}\left(b_{i}-\left\langle \boldsymbol{a}^{i},\boldsymbol{x}\right\rangle \right)^{2}}.\label{eq:Res_definition}
\end{equation}
\end{linenomath}
Note that each element of this problem set $\mathbb{S}$ specifies
an ordered sequence of $W$ blocks of linear equations of the form
$\left\langle \boldsymbol{a}^{i},\boldsymbol{x}\right\rangle =b_{i}$
where $\left\langle *,*\right\rangle $ denotes the inner product
in $\mathbb{R}^{J}$ (and thus $\mathbb{S}$ is an appropriate representation
of the so-called ``ordered subsets'' approach to tomographic reconstruction
\cite{HUDS94a}, as well as of other earlier-published block-iterative
methods that proposed essentially the same idea \cite{ELFV80a,EGGE81a,AHAR89a}).
The proximity function $Res$ on $\mathbb{S}$ is the \textit{residual}
that we get when a particular $\boldsymbol{x}$ is substituted into
all the equations of a particular problem $S$.

\subsection{Algorithms and outputs\label{sub:Algorithms-and-outputs}}

We now define the concept of an algorithm in the general context of
problem structures. For technical reasons that will become clear as
we proceed with our development, we introduce an additional set $\Delta$,
such that $\Omega\subseteq\Delta\subseteq\mathbb{R}^{J}$. (Both $\Omega$
and $\Delta$ are assumed to be known and fixed for any particular
problem structure $\left\langle \mathbb{T},\mathcal{P}r\right\rangle $.)
An \emph{algorithm }\textbf{$\mathbf{P}$ }for a problem structure
$\left\langle \mathbb{T},\mathcal{P}r\right\rangle $ assigns to each
problem $T\in\mathbb{T}$ an operator $\mathbf{P}_{T}:\Delta\rightarrow\Omega$.
This definition is used to define iterative processes that, for any
\textit{initial point} $\boldsymbol{x}\in\Omega,$ produce the (potentially)
infinite sequence $\left(\left(\mathbf{P}_{T}\right)^{k}\boldsymbol{x}\right)_{k=0}^{\infty}$
(that is, the sequence $\boldsymbol{x},\mathbf{P}_{T}\boldsymbol{x},\mathbf{P}_{T}\left(\mathbf{P}_{T}\boldsymbol{x}\right),\cdots$)
of points in $\Omega$. We discuss below how such a potentially infinite
process is terminated in practice.

Selecting $\Omega=\mathbb{R}_{+}^{J}$
and $\Delta=\mathbb{R}^{J}$ for the
problem structure $\left\langle \mathbb{S},Res\right\rangle $ of
the previous subsection, an example of an algorithm $\mathbf{R}$
is specified by
\begin{linenomath}
\begin{equation}
\mathbf{R}_{S}\boldsymbol{x}=\mathbf{QB}_{S_{W}}\cdots\mathbf{B}_{S_{1}}\boldsymbol{x},\label{eq:R_S_definition}
\end{equation}
\end{linenomath}
where $S$ is the problem specified above (\ref{eq:Res_definition})
and, for $1\leq w\leq W,$$\;$$\mathbf{B}_{S_{w}}:\Delta\rightarrow\Delta$
is defined by 
\begin{linenomath}
\begin{equation}
\mathbf{B}_{S_{w}}\boldsymbol{x}=\boldsymbol{x}+\frac{1}{\ell_{w}}\:\sum_{i=\ell_{1}+\ldots+\ell_{w-1}+1}^{\ell_{1}+\ldots+\ell_{w}}\frac{b_{i}-\left\langle \boldsymbol{a}^{i},\boldsymbol{x}\right\rangle }{\left\Vert \boldsymbol{a}^{i}\right\Vert ^{2}}\boldsymbol{a}^{i},\label{eq:P_s_w_definition}
\end{equation}
\end{linenomath}
where $\left\Vert \boldsymbol{a}\right\Vert $ denotes the norm of
the vector $\boldsymbol{a}$ in $\mathbb{R}^{J}$, and $\mathbf{Q}:\Delta\rightarrow\Omega$
is defined by 
\begin{linenomath}
\begin{equation}
\left(\mathbf{Q}\boldsymbol{x}\right)_{j}=\max\left\{ 0,\boldsymbol{x}_{j}\right\} ,\;\mbox{for }1\leq j\leq J.\label{eq:Q_definition}
\end{equation}
\end{linenomath}
Note that $\mathbf{R}_{S}:\Delta\rightarrow\Omega$. This specific
algorithm $\mathbf{R}$ is a typical example of the so-called block-iterative
methods mentioned above. Except for the presence of $\mathbf{Q}$
in (\ref{eq:R_S_definition}), which enforces nonnegativity of the
components, it is identical to an algorithm used and illustrated in
\cite{HERM08a}. With the $\mathbf{Q}$ absent from the definition
of the algorithm, $\Omega$ has to be the whole of $\mathbb{R}^{J}$;
the practical consequence of the presence versus the absence of $\mathbf{Q}$
in the tomographic application is illustrated in Subsection \ref{sub:Effects-of-variations}.
We note also that special cases of the presented algorithm include
the classical reconstruction methods ART (if $\ell_{w}=1,$ for $1\leq w\leq W$) and SIRT (if $W=1$); see, for example, Chapters
11 and 12 of \cite{HERM09a}.

For a problem structure $\left\langle \mathbb{T},\mathcal{P}r\right\rangle $,
a $T\in\mathbb{T}$, an $\varepsilon\in\mathbb{R}_{+}$ and a sequence
$R=\left(\boldsymbol{x}^{k}\right)_{k=0}^{\infty}$ of points in $\Omega$,
we use $O\left(T,\varepsilon,R\right)$ to denote the $\boldsymbol{x}\in\Omega$
that has the following properties: $\mathcal{P}r_{T}(\boldsymbol{x})\leq\varepsilon$\emph{
}and there is a nonnegative integer $K$ such that $\boldsymbol{x}^{K}=\boldsymbol{x}$
and, for all nonnegative integers $k<K$, $\mathcal{P}r_{T}\left(\boldsymbol{x}^{k}\right)>\varepsilon$.
Clearly, if there is such an $\boldsymbol{x}$, then it is unique.
If there is no such $\boldsymbol{x}$, then we say that $O\left(T,\varepsilon,R\right)$
is \textit{undefined,} otherwise we say that it is \textit{defined}.
The intuition behind this definition is the following: if we think
of $R$ as the (infinite) sequence of points that is produced by an
algorithm (intended for the problem $T$) without a termination criterion,
then $O\left(T,\varepsilon,R\right)$ is the \textit{output} produced
by that algorithm when we add to it instructions that make it terminate
as soon as it reaches a point that is $\varepsilon$-compatible with
$T$.

\subsection{Bounded perturbation resilience}

The notion of a \emph{bounded perturbations resilient }algorithm $\mathbf{P}$
for a problem structure $\left\langle \mathbb{T},\mathcal{P}r\right\rangle $
has been defined in a mathematically precise manner \cite{CENS10a}.
However, that definition is not satisfactory from the point of view
of applications in medical physics (or indeed in any area involving
noisy data), because it is useful only for problems $T$ for which
there is a perfectly compatible solution (that is, an $\boldsymbol{x}$
such that $\mathcal{P}r_{T}\left(\boldsymbol{x}\right)=0$). We therefore
extend here that notion as follows. An algorithm $\mathbf{P}$ for
a problem structure $\left\langle \mathbb{T},\mathcal{P}r\right\rangle $
is said to be \textit{strongly perturbation resilient} if, for all
$T\in\mathbb{T}$, 
\begin{enumerate}
\item there exists an $\varepsilon\in\mathbb{R}_{+}$ such that $O\left(T,\varepsilon,\left(\left(\mathbf{P}_{T}\right)^{k}\boldsymbol{x}\right)_{k=0}^{\infty}\right)$
is defined for every $\boldsymbol{x}\in\Omega$; 
\item for all $\varepsilon\in\mathbb{R}_{+}$ such that $O\left(T,\varepsilon,\left(\left(\mathbf{P}_{T}\right)^{k}\boldsymbol{x}\right)_{k=0}^{\infty}\right)$
is defined for every $\boldsymbol{x}\in\Omega$, we also have that
$O\left(T,\varepsilon',R\right)$ is defined for every $\varepsilon'>\varepsilon$
and for every sequence $R=\left(\boldsymbol{x}^{k}\right)_{k=0}^{\infty}$
of points in $\Omega$ generated by 
\begin{linenomath}
\begin{equation}
\boldsymbol{x}^{k+1}=\mathbf{P}_{T}\left(\boldsymbol{x}^{k}+\beta_{k}\boldsymbol{v}^{k}\right),\:\mathrm{for\: all\:}k\geq0,\label{eq:perturbations}
\end{equation}
\end{linenomath}
where $\beta_{k}\boldsymbol{v}^{k}$ are \textit{bounded perturbations},
meaning that the sequence $\left(\beta_{k}\right)_{k=0}^{\infty}$
of nonnegative real numbers is \textit{summable} (that is, ${\displaystyle \sum\limits _{k=0}^{\infty}}\beta_{k}\,<\infty$),
the sequence $\left(\boldsymbol{v}^{k}\right)_{k=0}^{\infty}$ of
vectors in $\mathbb{R}^{J}$ is bounded and, for all $k\geq0$, $\boldsymbol{x}^{k}+\beta_{k}\boldsymbol{v}^{k}\in\Delta$. 
\end{enumerate}
$\quad\;$In less formal terms, the second of these properties says
that for a strongly perturbation resilient algorithm we have that,
for every problem and any nonnegative real number $\varepsilon$,
if it is the case that for all initial points from $\Omega$ the infinite
sequence produced by the algorithm contains an $\varepsilon$-compatible
point, then it will also be the case that all perturbed sequences
satisfying (\ref{eq:perturbations}) contain an $\varepsilon'$-compatible
point, for any $\varepsilon'>\varepsilon$.

Having defined the notion of a strongly perturbation resilient algorithm,
we next show that this notion is of relevance to problems in medical
physics. We illustrate the use of this in tomography in the next section.
We first need to introduce some mathematical concepts.

Given an algorithm $\mathbf{P}$ for a problem structure $\left\langle \mathbb{T},\mathcal{P}r\right\rangle $
and a $T\in\mathbb{T}$, we say that $\mathbf{P}$ is \textit{convergent
for $T$} if, for every $\boldsymbol{x}\in\Omega$, there exists a
unique $\boldsymbol{y}\left(\boldsymbol{x}\right)\in\Omega$ such
that, $lim_{k\rightarrow\infty}\left(\mathbf{P}_{T}\right)^{k}\boldsymbol{x}=\boldsymbol{y}\left(\boldsymbol{x}\right)$,
meaning that for every positive real number $\delta$, there exist
a nonnegative integer $K$, such that $\left\Vert \left(\mathbf{P}_{T}\right)^{k}\boldsymbol{x}-\boldsymbol{y}\left(\boldsymbol{x}\right)\right\Vert \leq\delta$,
for all nonnegative integers $k\geq K$. If, in addition, there exists
a $\gamma\in\mathbb{R}_{+}$ such that $\mathcal{P}r_{T}\left(\boldsymbol{y}\left(\boldsymbol{x}\right)\right)\leq\gamma$,
for every $\boldsymbol{x}\in\Omega$, then we say that $\mathbf{P}$
is \textit{boundedly convergent for $T$}.

A function $f:\Omega\rightarrow\mathbb{R}$ is \textit{uniformly continuous}
if, for every $\varepsilon>0$ there exists a $\delta>0$, such that,
for all $\boldsymbol{x},\boldsymbol{y}\in\Omega$, $\left|f(\boldsymbol{x})-f(\boldsymbol{y})\right|\leq\varepsilon$
provided that $\left\Vert \boldsymbol{x}-\boldsymbol{y}\right\Vert \leq\delta$.
An example of a uniformly continuous function is $Res_{S}$ of (\ref{eq:Res_definition}),
for any $S\in\mathbb{S}$. This can be proved by observing that the
right-hand side of (\ref{eq:Res_definition}) can be rewritten in
vector/matrix form as $\left\Vert \boldsymbol{b-Ax}\right\Vert $
and then selecting, for any given $\varepsilon>0$, $\delta$ to be
$\varepsilon/\left\Vert \boldsymbol{A}\right\Vert $, where $\left\Vert \boldsymbol{A}\right\Vert $
denotes the matrix norm of $\boldsymbol{A}$.

An operator ${\bf O}:\Delta\rightarrow\Omega$, is \emph{nonexpansive}
if $\left\Vert {\bf O}\boldsymbol{x}-{\bf O}\boldsymbol{y}\right\Vert \leq\left\Vert \boldsymbol{x}-\boldsymbol{y}\right\Vert $,
for all $\boldsymbol{x},\boldsymbol{y}\in\Delta$. An example of a
nonexpansive operator is the $\mathbf{R}_{S}$ of (\ref{eq:R_S_definition}).
The proof of this is also simple. It follows from discussions regarding
similar claims in \cite{BUTN07a} that the $\mathbf{B}_{S_{w}}:\mathbb{R}^{J}\rightarrow\mathbb{R}^{J}$
of (\ref{eq:P_s_w_definition}) is a nonexpansive operator, for $1\leq w\leq W,$
and that the operator $\mathbf{Q}$ of (\ref{eq:Q_definition}) is
also nonexpansive. Obviously, a sequential application of nonexpansive
operators results in a nonexpansive operator and thus $\mathbf{R}_{S}$
is nonexpansive.

Now we state an important new result that gives sufficient conditions
for strong perturbation resilience: \texttt{If $\mathbf{P}$ is an
algorithm for a problem structure $\left\langle \mathbb{T},\mathcal{P}r\right\rangle $
such that, for all $T\in\mathbb{T}$, $\mathbf{P}$ is boundedly convergent
for $T$, $\mathcal{P}r_{T}:\Omega\rightarrow\mathbb{R}$ is uniformly
continuous and $\mathbf{P}_{T}:\Delta\rightarrow\Omega$ is nonexpansive,
then $\mathbf{P}$ is strongly perturbation resilient.} The importance
of this result lies in the fact that the rather ordinary condition
of uniform continuity for the proximity function and the reasonable
conditions of bounded convergence and nonexpansiveness of the algorithmic
operators guarantee that we end up with a strongly perturbation resilient
algorithm. The proof of this new result involves some mathematical
technicalities and is therefore presented in the Appendix as Theorem
1.

\subsection{Optimization criterion and nonascending vector}

Now suppose, as is indeed the case for the constrained optimization
problems discussed in the previous section, that in addition to a
problem structure $\left\langle \mathbb{T},\mathcal{P}r\right\rangle $
we are also provided with an optimization criterion, which is specified
by a function $\phi:\Delta\rightarrow\mathbb{R}$, with the convention
that a point in $\Delta$ for which the value of $\phi$ is smaller
is considered \textit{superior} (from the point of view of our application)
to a point in $\Delta$ for which the value of $\phi$ is larger.
In the tomography context, any of the functions of $\boldsymbol{x}$
that are listed as a ``secondary optimization criterion'' (an alternative
name is a ``regularizer'') in Section 6.4 of\cite{HERM09a} is an
acceptable choice for the optimization criterion $\phi$. These include
weighted norms, the negative of Shannon's entropy and total variation.
It is the last of these that we discuss in detail in the illustrative
example below. The essential idea of the \emph{superiorization methodology
}presented in this paper is to make use of the perturbations of (\ref{eq:perturbations})
to transform a strongly perturbation resilient algorithm that seeks
a constraints-compatible solution into one whose outputs are equally
good from the point of view of constraints-compatibility, but are
superior according to the optimization criterion. We do this by producing
from the algorithm another one, called its \textit{superiorized} version,
by making sure not only that the $\beta_{k}\boldsymbol{v}^{k}$ are
bounded perturbations, but also that $\phi\left(\boldsymbol{x}^{k}+\beta_{k}\boldsymbol{v}^{k}\right)\leq\phi\left(\boldsymbol{x}^{k}\right)$,
for all $k\geq0$.

In order to ensure this we introduce a new concept (closely related
to the concept of a ``descent direction'' that is widely used in
optimization). Given a function $\phi:\Delta\rightarrow\mathbb{R}$
and a point $\boldsymbol{x}\in\Delta$, we say that a vector $\boldsymbol{d}\in\mathbb{R}^{J}$
is \textit{nonascending} for $\phi$ at $\boldsymbol{x}$ if $\left\Vert \boldsymbol{d}\right\Vert \leq1$
and 
\begin{linenomath}
\begin{equation}
\begin{array}{r}
\mathrm{there\: is\: a\:\delta>0\: such\: that\: for\: all\:}\lambda\in\left[0,\delta\right],\\
\left(\boldsymbol{x}+\lambda\boldsymbol{d}\right)\in\Delta\:\mathrm{and}\:\phi\left(\boldsymbol{x}+\lambda\boldsymbol{d}\right)\leq\phi\left(\boldsymbol{x}\right).
\end{array}\label{eq:nonascending}
\end{equation}
\end{linenomath}

Note that irrespective of the choices of $\phi$ and $\boldsymbol{x}$,
there is always at least one nonascending vector $\boldsymbol{d}$
for $\phi$ at $\boldsymbol{x}$, namely the zero-vector, all of whose
components are zero. This is a useful fact for proving results concerning
the guaranteed behavior of our proposed procedures. However, in order
to steer our algorithms toward a point at which the value of $\phi$
is small, we need to find a $\boldsymbol{d}$ such that $\phi\left(\boldsymbol{x}+\lambda\boldsymbol{d}\right)<\phi\left(\boldsymbol{x}\right)$
rather than just $\phi\left(\boldsymbol{x}+\lambda\boldsymbol{d}\right)\leq\phi\left(\boldsymbol{x}\right)$
as in (\ref{eq:nonascending}). In some earlier papers on superiorization
\cite{CENS10a,DAVI09a,HERM08a,NIKA12a,BUTN07a} it was assumed that
$\Delta=\mathbb{R}^{J}$ and that $\phi$ is a convex function. This
implied that, for any point $\boldsymbol{x}\in\Delta$, $\phi$ had
a subgradient $\boldsymbol{g}\in\mathbb{R}^{J}$ at the point $\boldsymbol{x}$.
It was suggested that if there is such a $\boldsymbol{g}$ with a
positive norm, then $\boldsymbol{d}$ should be chosen to be $-\boldsymbol{g}/\left\Vert \boldsymbol{g}\right\Vert $,
otherwise $\boldsymbol{d}$ should be chosen to be the zero vector.
However, there are approaches (not involving subgradients) to selecting
an appropriate $\boldsymbol{d}$; an example can be found in \cite{GARD11a}
in which $\boldsymbol{d}$ is found without using subgradients for
the case when $\phi$ is the $\ell_{1}$-norm of the Haar transform.
The method we used for selecting a nonascending vector in the experiments
reported in this paper is specified at the end of Subsection \ref{sub:Application-to-tomography}.

\subsection{\label{sub:Superiorized-version-of}Superiorized version of an algorithm}

We now make precise the ingredients needed for transforming an algorithm
into its superiorized version. Let $\Omega$ and $\Delta$ be the
underlying sets for a problem structure $\left\langle \mathbb{T},\mathcal{P}r\right\rangle $\textbf{
(}$\Omega\subseteq\Delta\subseteq\mathbb{R}^{J}$, as discussed at
the beginning of Subsection \ref{sub:Algorithms-and-outputs}), $\mathbf{P}$\textbf{
}be an algorithm for $\left\langle \mathbb{T},\mathcal{P}r\right\rangle $
and $\phi:\Delta\rightarrow\mathbb{R}$. The following description
of the Superiorized Version of Algorithm $\mathbf{P}$ produces, for
any problem $T\in\mathbb{T}$, a sequence $R_{T}=\left(\boldsymbol{x}^{k}\right)_{k=0}^{\infty}$
of points in $\Omega$ for which, for all $k\geq0$, (\ref{eq:perturbations})
is satisfied. We show this to be true, for any algorithm $\mathbf{P},$
after the description of the Superiorized Version of Algorithm $\mathbf{P}$.
Furthermore, since the sequence $R_{T}$ is steered by Superiorized
Version of Algorithm $\mathbf{P}$ toward a reduced value of $\phi$,
there is an intuitive expectation that the output of the superiorized
version is likely to be superior (from the point of view of the optimization
criterion $\phi$) to the output of the original unperturbed algorithm.
This last statement is not precise and so it cannot be proved in a
mathematical sense for an arbitrary algorithm $\mathbf{P}$; however,
that should not stop us from applying the easy procedure given below
for automatically producing the Superiorized Version of $\mathbf{P}$
and experimentally checking whether it indeed provides us with outputs
superior to those of the original algorithm. The well-demonstrated
nature of heuristic optimization approaches is that they often work
in practice even when their performance cannot be guaranteed to be
optimal\cite{RARD01a,WERN00a,ZANA81a}.

Nevertheless, we can push our theory further than the hope expressed
in the last paragraph, by considering superiorized versions of algorithms
that satisfy some condition. In this paper, the condition that we
discuss is strong perturbation resilience. We show below that if $\mathbf{P}$
is strongly perturbation resilient, then, for any problem $T\in\mathbb{T}$,
a sequence $R_{T}$ produced by its superiorized version has the following
desirable property: For all $\varepsilon\in\mathbb{R}_{+}$, if $O\left(T,\varepsilon,\left(\left(\mathbf{P}_{T}\right)^{k}\boldsymbol{x}\right)_{k=0}^{\infty}\right)$
is defined for every $\boldsymbol{x}\in\Omega$, then $O\left(T,\varepsilon',R_{T}\right)$
is also defined for every $\varepsilon'>\varepsilon$; in other words,
the Superiorized Version of Algorithm $\mathbf{P}$ provides an $\varepsilon'$-compatible
output. As stated above, the advantage of the superiorized version
is that its output is likely to be superior to the output of the original
unperturbed algorithm. We point out that strong perturbation resilience
is a sufficient, but not necessary, condition for guaranteeing such
desirable behavior of the superiorized version, finding additional
sufficient conditions and proving that algorithms that we wish to
superiorize satisfy such conditions is part of our ongoing research.

The superiorized version assumes that we have available a summable
sequence $\left(\gamma_{\ell}\right)_{\ell=0}^{\infty}$ of positive
real numbers (for example, $\gamma_{\ell}=a^{\ell}$, where $0<a<1$)
and it generates, simultaneously with the sequence $\left(\boldsymbol{x}^{k}\right)_{k=0}^{\infty}$,
sequences $\left(\boldsymbol{v}^{k}\right)_{k=0}^{\infty}$ and $\left(\beta_{k}\right)_{k=0}^{\infty}$.
The latter is generated as a subsequence of $\left(\gamma_{\ell}\right)_{\ell=0}^{\infty}$,
resulting in a summable sequence $\left(\beta_{k}\right)_{k=0}^{\infty}$.
The algorithm further depends on a specified initial point $\boldsymbol{\bar{x}}\in\Omega$
and on a positive integer $N$. It makes use of a logical variable
called \emph{loop.}

\textbf{\textit{\emph{Superiorized Version of Algorithm $\mathbf{P}$}}} 
\begin{enumerate}
\item \textbf{set} $k=0$ 
\item \textbf{set} $\boldsymbol{x}^{k}=\boldsymbol{\bar{x}}$ 
\item \textbf{set} $\ell=-1$ 
\item \textbf{repeat} 
\item $\qquad$\textbf{set} $n=0$ 
\item $\qquad$\textbf{set} $\boldsymbol{x}^{k,n}=\boldsymbol{x}^{k}$ 
\item $\qquad$\textbf{while $n<N$} 
\item $\qquad\qquad$\textbf{set $\boldsymbol{v}^{k,n}$ }to be a nonascending
vector for $\phi$ at $\boldsymbol{x}^{k,n}$ 
\item $\qquad\qquad$\textbf{set} \emph{loop=true} 
\item $\qquad\qquad$\textbf{while}\emph{ loop} 
\item $\qquad\qquad\qquad$\textbf{set $\ell=\ell+1$} 
\item $\qquad\qquad\qquad$\textbf{set} $\beta_{k,n}=\gamma_{\ell}$ 
\item $\qquad\qquad\qquad$\textbf{set} $\boldsymbol{z}=\boldsymbol{x}^{k,n}+\beta_{k,n}\boldsymbol{v}^{k,n}$ 
\item $\qquad\qquad\qquad$\textbf{if $\boldsymbol{z}\in\Delta$ $\mathbf{\: and\:}$
$\phi\left(\boldsymbol{z}\right)\leq\phi\left(\boldsymbol{x}^{k}\right)$
$\mathbf{then}$}
\item $\qquad\qquad\qquad\qquad$\textbf{set $n=n+1$} 
\item $\qquad\qquad\qquad\qquad$\textbf{set $\boldsymbol{x}^{k,n}=\boldsymbol{z}$} 
\item $\qquad\qquad\qquad\qquad$\textbf{set }\emph{loop = false} 
\item $\qquad$\textbf{set $\boldsymbol{x}^{k+1}=\mathbf{P}_{T}\boldsymbol{x}^{k,N}$} 
\item $\qquad$\textbf{set $k=k+1$} 
\end{enumerate}
Next we analyze the behavior of the Superiorized Version of Algorithm
$\mathbf{P}$.

The iteration number $k$ is set to 0 in (i) and $\boldsymbol{x}^{k}=\boldsymbol{x}^{0}$$\:$
is set to its initial value $\bar{\boldsymbol{x}}$ in (ii). The integer
index $\ell$ for picking the next element from the sequence $\left(\gamma_{\ell}\right)_{\ell=0}^{\infty}$
is initialized to $-1$ by line (iii), it is repeatedly increased
by line (xi). The lines (v) - (xix) that follow the $\mathbf{repeat}$
in (iv) perform a complete iterative step from $\boldsymbol{x}^{k}$
to $\boldsymbol{x}^{k+1}$, infinite repetitions of such steps provide
the sequence $R_{T}=\left(\boldsymbol{x}^{k}\right)_{k=0}^{\infty}$.
During one iterative step, there is one application of the operator
$\mathbf{P}_{T}$, in line (xviii), but there are $N$ steering steps
aimed at reducing the value of $\phi$; the latter are done by lines
(v) - (xvii). \texttt{These lines produce a sequence of points $\boldsymbol{x}^{k,n}$,
where $0\leq n\leq N$ with $\boldsymbol{x}^{k,0}=\boldsymbol{x}^{k}$,
$\boldsymbol{x}^{k,n}\in\Delta$ and $\phi\left(\boldsymbol{x}^{k,n}\right)\leq\phi\left(\boldsymbol{x}^{k}\right)$.}

We prove the truth of the last sentence by induction on the nonnegative
integers. For $n=0$, we have by lines (v) and (vi) that $\boldsymbol{x}^{k,0}=\boldsymbol{x}^{k}$.
But $\boldsymbol{x}^{k}\in\Omega$ , since it is either $\boldsymbol{\bar{x}}$
that is assumed to be in $\Omega$ due to lines (i) and (ii) or it
is in the range $\Omega$ of $\mathbf{P}_{T}$ due to lines (xviii)
and (xix). Now we assume, for any $0\leq n<N$, that $\boldsymbol{x}^{k,n}\in\Delta$
and $\phi\left(\boldsymbol{x}^{k,n}\right)\leq\phi\left(\boldsymbol{x}^{k}\right)$
and show that lines (viii) - (xvii) perform a computation that leads
from $\boldsymbol{x}^{k,n}$ to an $\boldsymbol{x}^{k,n+1}\in\Delta$
that satisfies $\phi\left(\boldsymbol{x}^{k,n+1}\right)\leq\phi\left(\boldsymbol{x}^{k}\right)$.
To see this, observe that line (viii) sets $\boldsymbol{v}^{k,n}$
to be a nonascending vector for $\phi$ at $\boldsymbol{x}^{k,n}$,
which implies that (\ref{eq:nonascending}) is satisfied with $\boldsymbol{x}=\boldsymbol{x}^{k,n}$
and $\boldsymbol{d}=\boldsymbol{v}^{k,n}$. Line (ix) sets \textit{loop}
to \textit{true}, and it remains \textit{true} while searching for
the desired $\boldsymbol{x}^{k,n+1}$, by repeatedly executing the
loop sequence that follows line (x). In this sequence, line (xi) increases
$\ell$ by 1 and line (xii) sets $\beta_{k,n}$ to $\gamma_{\ell}$.
Thus for the vector $\boldsymbol{z}$ defined by line (xiii), $\boldsymbol{z}\in\Delta$
and $\phi\left(\boldsymbol{z}\right)\leq\phi\left(\boldsymbol{x}^{k,n}\right)$,
provided that $\beta_{k,n}$ is not greater than the $\delta$ in
(\ref{eq:nonascending}). Since $\left(\gamma_{\ell}\right)_{\ell=0}^{\infty}$
is a summable sequence of positive real numbers, there must be a positive
integer $L$ such that $\gamma_{\ell}\leq\delta$, for all $\ell\geq L$.
This implies that if we applied lines (xi) - (xiii) often enough,
we would reach a vector $\boldsymbol{z}$ that satisfies $\boldsymbol{z}\in\Delta$
and $\phi\left(\boldsymbol{z}\right)\leq\phi\left(\boldsymbol{x}^{k,n}\right)$.
If the condition in line (xiv) is not satisfied when the process gets
to it, then lines (xi) - (xiii) are again executed and eventually
we get a vector $\boldsymbol{z}$ for which the condition in line
(xiv) is satisfied due to the induction hypothesis that $\phi\left(\boldsymbol{x}^{k,n}\right)\leq\phi\left(\boldsymbol{x}^{k}\right)$.
By lines (xv) and (xvi) we see that at that time $\boldsymbol{x}^{k,n+1}$
is set to $\boldsymbol{z}$ and so we obtain that $\boldsymbol{x}^{k,n+1}\in\Delta$
and $\phi\left(\boldsymbol{x}^{k,n+1}\right)\leq\phi\left(\boldsymbol{x}^{k}\right)$,
as desired. Line (xvii) sets \textit{loop} to \textit{false} and so
control is returned to line (vii). When this happens for the $N$th
time, it will be the case that $n=N$ and therefore line (xviii) is
used to produce $\boldsymbol{x}^{k+1}\in\Omega$ and the increasing
of $k$ by line (xix) allows us then to move on to the next iterative
step. Infinite repetition of such steps produces the sequence $R_{T}=\left(\boldsymbol{x}^{k}\right)_{k=0}^{\infty}$
of points in $\Omega$.

We now show that if $O\left(T,\varepsilon,\left(\left(\mathbf{P}_{T}\right)^{k}\boldsymbol{x}\right)_{k=0}^{\infty}\right)$
is defined for every $\boldsymbol{x}\in\Omega$, then, for any $\varepsilon'>\varepsilon$,
the Superiorized Version of Algorithm $\mathbf{P}$ produces an $\varepsilon'$-compatible
output. Since $\mathbf{P}$ is assumed to be strongly perturbation
resilient, this desired result follows if we can show that there exists
a summable sequence $\left(\beta_{k}\right)_{k=0}^{\infty}$ of nonnegative
real numbers and a bounded sequence $\left(\boldsymbol{v}^{k}\right)_{k=0}^{\infty}$
of vectors in $\mathbb{R}^{J}$ such that (\ref{eq:perturbations})
is satisfied for all $k\geq0.$ In view of line (xviii), this is achieved
if we can define the $\beta_{k}$ and the $\boldsymbol{v}^{k}$ so
that $\boldsymbol{x}^{k,N}=\boldsymbol{x}^{k}+\beta_{k}\boldsymbol{v}^{k}$.
This is done by setting 
\begin{linenomath}
\begin{equation}
\beta_{k}=\max\left\{ \beta_{k,n}\,|\,0\leq n<N\right\} ,\label{eq:beta_k}
\end{equation}
\end{linenomath}
 
\begin{equation}
\boldsymbol{v}^{k}=\sum_{n=0}^{N-1}\frac{\beta_{k,n}}{\beta_{k}}\boldsymbol{v}^{k,n}.\label{eq:v^k}
\end{equation}
That these assignments result in $\boldsymbol{x}^{k,N}=\boldsymbol{x}^{k}+\beta_{k}\boldsymbol{v}^{k}$
follows from lines (v) - (xvii). From line (xii) follows that $\left(\beta_{k}\right)_{k=0}^{\infty}$
is a subsequence of $\left(\gamma_{\ell}\right)_{\ell=0}^{\infty}$
and, hence, it is a summable sequence of nonnegative real numbers.
Since each $\left\Vert \boldsymbol{v}^{k,n}\right\Vert \leq1$ by
the definition of a nonascending vector, it follows from (\ref{eq:beta_k})
and (\ref{eq:v^k}) that $\left\Vert \boldsymbol{v}^{k}\right\Vert \leq N$
and so $\left(\boldsymbol{v}^{k}\right)_{k=0}^{\infty}$ is bounded.
Part of the condition expressed in (\ref{eq:perturbations}) is that,
for all $k\geq0$, $\boldsymbol{x}^{k}+\beta_{k}\boldsymbol{v}^{k}\in\Delta$.
This follows from the fact that $\boldsymbol{x}^{k,N}=\boldsymbol{x}^{k}+\beta_{k}\boldsymbol{v}^{k}$
is assigned its value by line (xvi), but only if the condition expressed
in line (xiv) is satisfied.

In conclusion, we have shown that the superiorized version of a strongly
perturbation resilient algorithm produces outputs that are essentially
as constraints-compatible as those produced by the original version
of the algorithm. However, due to the repeated steering of the process
by lines (vii) - (xvii) toward reducing the value of the optimization
criterion $\phi$, we can expect that the output of the superiorized
version will be superior (from the point of view of $\phi$) to the
output of the original algorithm.

\subsection{\label{sub:Information-on-performance}Information on performance
comparison with MAP methods}

Using our notation, the constrained minimization formulation that
we are considering is: Given an $\varepsilon\in\mathbb{R}_{+}$,
\begin{linenomath}
\begin{equation}
\textrm{minimize }\phi(\boldsymbol{x})\textrm{, subject to }\mathcal{P}r_{T}\left(\boldsymbol{x}\right)\leq\varepsilon.\label{eq:constrained_minimization}
\end{equation}
\end{linenomath}
The aim of superiorization is not identical with the aim of constrained
minimization in (\ref{eq:constrained_minimization}). One difference
is that $\varepsilon$ is not ``given'' in the superiorization context.
The superiorization of an algorithm produces a sequence and, for any
$\varepsilon$, the associated output of the algorithm is considered
to be the first $\boldsymbol{x}$ in the sequence for which $\mathcal{P}r_{T}\left(\boldsymbol{x}\right)\leq\varepsilon$.
The other difference is that we do not claim that this output is a
minimizer of $\phi$ among all points that satisfy the constraint,
but hope only that it is usually an $\boldsymbol{x}$ for which $\phi(\boldsymbol{x})$
is at the small end of its range of values over the set of constraint-satisfying
points. This latter difference is generally shared by comparisons
of a heuristic approach with an exact approach to solving a constrained
minimization problem.

The MAP (or regularized) formulation of a physical problem that leads
to the constrained minimization problem (\ref{eq:constrained_minimization})
is the unconstrained minimization problem of the form: Given a $\beta\in\mathbb{R}_{+}$,
\begin{linenomath}
\begin{equation}
\textrm{minimize }\left[\phi(\boldsymbol{x})+\beta\mathcal{P}r_{T}\left(\boldsymbol{x}\right)\right].\label{eq:MAP_formulation}
\end{equation}
\end{linenomath}
Formulations of both kinds (i.e, the ones of (\ref{eq:constrained_minimization})
and of (\ref{eq:MAP_formulation})) are widely used for solving medical
physics problems and the question ``Which
of these two formulations leads to faster or better solutions of the
underlying physical problem?'' is open. Examples of both formulations
with various choices for $\mathcal{P}r_{T}$ and $\phi$ are listed
in the beginning parts of the paper of Goldstein and Osher\cite{GOLD09a}.

We now return to the question raised near the end of Section \ref{sec:Introduction}:
Will superiorization produce superior results
to those produced by contemporary MAP methods or is it faster than
the better of such methods? As yet, there is very little information
available regarding this general question; in fact, we are aware of
only one published study\cite{CENS12a}. That study compared a superiorization
algorithm with the algorithm of Goldstein and Osher that they refer
to as TwIST \cite{BIOU07a} with split Bregman \cite{GOLD09a} as
the substep, which is indeed a contemporary method that uses the MAP
formulation. (For example, see the discussion of the split Bregman
method in \cite{ABAS11a}.) The problem $S$ to which the two algorithms
were applied was one from the tomographic problem set $\mathbb{S}$
defined in (\ref{eq:S_definition}). $Res_{S}$ as defined in (\ref{eq:Res_definition})
was used as the proximity function and total variation, $TV$ as defined
below in (\ref{eq:TV}), was the choice for $\phi$. It is reported
in \cite{CENS12a} that for the outputs of the two algorithms that
were being compared, the values of $Res_{S}$
and $TV$ were very similar, but the superiorization algorithm produced
its output four times faster than the MAP method.

\section{\label{sec:An-Illustrative-ExampleResults}An Illustrative Example}

\subsection{Application to tomography\label{sub:Application-to-tomography}}

We use \textit{tomography} to refer to the process of reconstructing
a function over a Euclidean space from estimated values of its integrals
along lines (that are usually, but not necessarily, straight). The
particular reconstruction processes to which our discussion applies
are the \textit{series expansion methods}, see Section 6.3 of \cite{HERM09a},
in which it is assumed that the function to be reconstructed can be
approximated by a linear combination of a finite number (say $J$)
of basis functions and the reconstruction task becomes one of estimating
the coefficients of the basis functions in the expansion. Sometimes,
prior knowledge about the nature of the function to be reconstructed
allows us to confine the sought-after vector $\boldsymbol{x}$ of
coefficients to a subset $\Omega$ of $\mathbb{R}^{J}$ (such as the
nonnegative orthant $\mathbb{R}_{+}^{J}$). We use $i$ to index the
lines along which we integrate, $\boldsymbol{a}^{i}\in\mathbb{R}^{J}$
to denote the vector whose $j$th component is the integral of the
$j$th basis function along the $i$th line, and $b_{i}$ to denote
the measured integral of the function to be reconstructed along the
$i$th line. Under these circumstances the constraints come from the
desire that, for each of the lines, $\left\langle \boldsymbol{a}^{i},\boldsymbol{x}\right\rangle $
should be close (in some sense) to $b_{i}$.

To make this concrete, consider (\ref{eq:S_definition}). Such a description
of the constraints arises in tomography by grouping the lines of integration
into $W$ blocks, with $\ell_{w}$ lines in the $w$th block. Such
groupings often (but not always) are done according to some geometrical
condition on the lines (for example, in case of straight lines, we
may decide that all the lines that are parallel to each other form
one block). In this framework the proximity function $Res$ defined
by (\ref{eq:Res_definition}) provides a reasonable measure of the
incompatibility of a vector $\boldsymbol{x}$ with the constraints.
The algorithm $\mathbf{R}$ described by (\ref{eq:R_S_definition})
- (\ref{eq:Q_definition}) is applicable to this concrete formulation.

There are many optimization criteria that have been used in tomography,
see Section 6.4 of \cite{HERM09a}, here we discuss the one called
\textit{total variation} ($TV$), whose use has been popular in medical
physics recently, see as examples \cite{KIM11a,PENF10a,SIDK11a,COMB04a,BIAN10A,DEFR11a,SIDK08a}.
The definition of $TV$ that we use here requires a certain way of
selecting the basis functions. It is assumed that the function to
be reconstructed is defined in the plane $\mathbb{R}^{2}$ and is
zero-valued outside a square-shaped region in the plane. This region
is subdivided into $J$ smaller equal-sized squares (\textit{pixels})
and the $J$ basis functions are defined by having value one in exactly
one pixel and value zero everywhere else. We index the pixels by $j$
and we let $C$ denote the set of all indices of pixels that are not
in the rightmost column or the bottom row of the pixel array. For
any pixel with index $j$ in $C$, let $r(j)$ and $b(j)$ be the
index of the pixel to its right and below it, respectively. We define
$TV:\mathbb{R}^{J}\rightarrow\mathbb{R}$ by 
\begin{linenomath}
\begin{equation}
TV(\boldsymbol{x})=\sum_{j\in C}\sqrt{\left(x_{j}-x_{r(j)}\right)^{2}+\left(x_{j}-x_{b(j)}\right)^{2}}.\label{eq:TV}
\end{equation}
\end{linenomath}

The method we adopted to generate a nonascending vector for the $TV$
function at an $\boldsymbol{x}\in\mathbb{R}^{J}$ is based on Theorem
2 of the Appendix. It is applicable since $TV:\mathbb{R}^{J}\rightarrow\mathbb{R}$
is a convex function; see, for example, the end of the Proof of Proposition
1 of \cite{COMB04a}. Now consider an integer $j'$ such that $1\leq j'\leq J$.
Looking at the sum in (\ref{eq:TV}), we see that $x_{j'}$ appears
in at most three terms, in which $j'$ must be either $j$, or $r(j),$
or $b(j)$ for some $j\in C.$ By taking the formal partial derivatives
of these three terms, we see that $\frac{\partial TV}{\partial x_{j'}}(\boldsymbol{x})$
is well-defined if the denominator in the formal derivative of any
of the three terms is not zero for $\boldsymbol{x}$. In view of this,
we define the $\boldsymbol{g}$ in Theorem 2 as follows. If the denominator
in any of the three formal partial derivatives with respect to $x_{j'}$
has an absolute value less than a very small positive number (we used
$10^{-20}$ ), then we set $g_{j'}$ to zero, otherwise we set it
to $\frac{\partial TV}{\partial x_{j'}}(\boldsymbol{x})$. Clearly
the resulting $\boldsymbol{g}\in\mathbb{R}^{J}$ satisfies the condition
in Theorem 2 and hence provides a $\boldsymbol{d}$ that is a nonascending
vector for $TV$ at $\boldsymbol{x}$.

Previously reported reconstructions using $TV$-superiorization selected
the $\boldsymbol{d}$ using subgradients as discussed in the paragraph
following (\ref{eq:nonascending}); such a $\boldsymbol{d}$ is not
guaranteed to be a nonascending vector for the $TV$ function. What
we are proposing here is not only mathematically rigorous (in the
sense that it is guaranteed to produce a nonascending vector for the
$TV$ function), but it can also lead to a better reconstructions,
as illustrated in Subsection \ref{sub:Effects-of-variations}.

\subsection{\label{sub:The-data-sets}The data generation for the experiments}

The data sets used in the experiments reported in this paper were
generated in such a way that they share the noise-characteristics
of CT scanners when used for scanning the human head and brain; as
discussed, for example, in Chapter 5 of \cite{HERM09a}. They were
generated using the software SNARK09 \cite{DAVI09b}.

\begin{figure}
\includegraphics[scale=0.5]{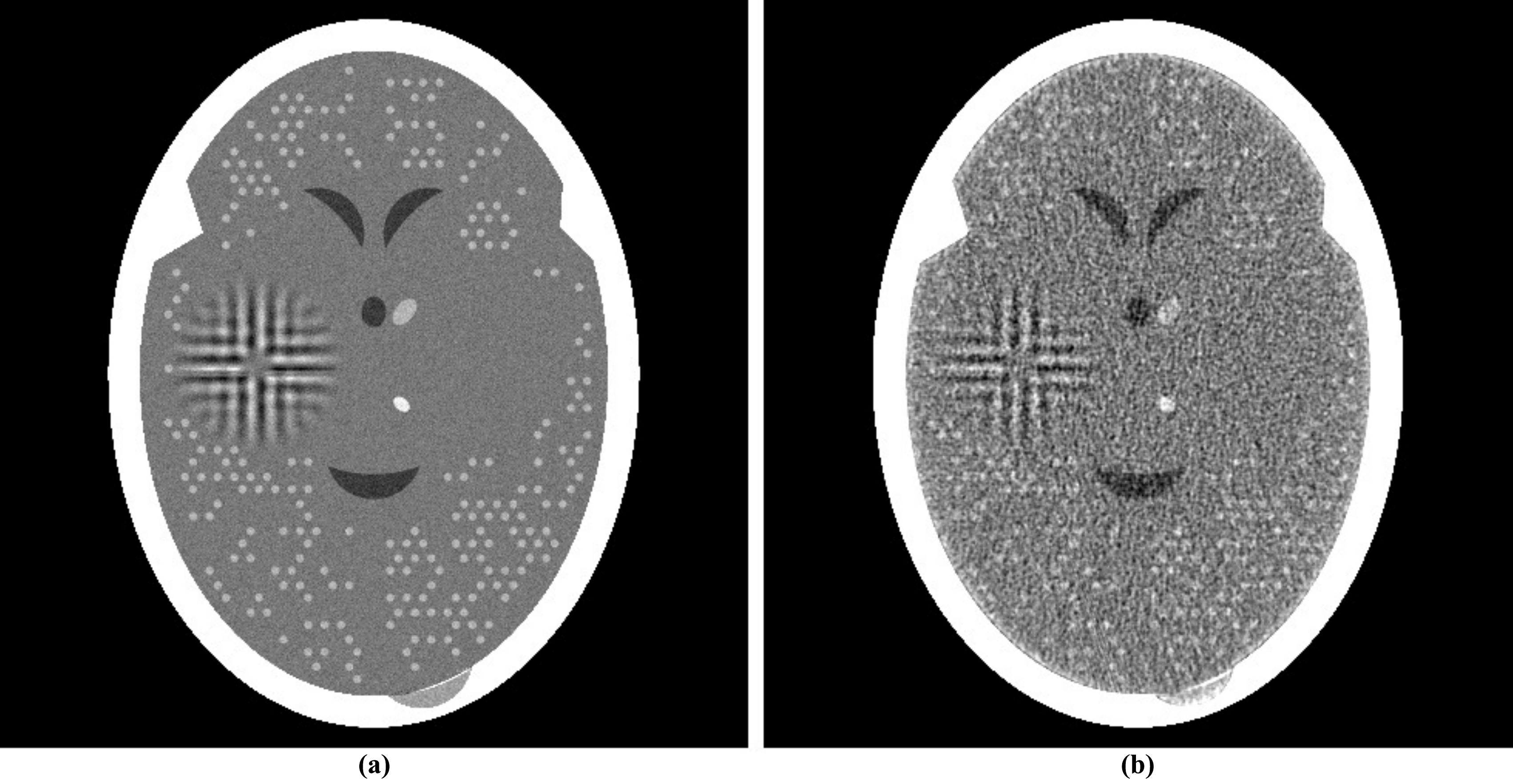}\caption{\label{fig:(a)-A-head}(a) A head phantom. (b) Reconstruction of the
head phantom from realistically simulated projection data for 360
views using ART with blob basis functions.}
\end{figure}

The head phantom that was used for data generation is based on an
actual cross-section of the human head. It is described as a collection
of geometrical objects (such as ellipses, triangles and segments of
circles) whose combination accurately resembles the anatomical features
of the actual head cross-section. In addition, the basic phantom contains
a large tumor. The actual phantom used was obtained by a random variation
of the basic phantom, by incorporating into it local inhomogeneities
and small low-contrast tumors at random locations. This phantom is
represented by the image in figure \ref{fig:(a)-A-head}. That image
comprises $485\times485$ pixels each of size 0.376 mm by 0.376 mm.
The values assigned to the pixels are obtained by an $11\times11$
sub-sampling of the pixels and averaging the values assigned to the
sub-samples by the geometrical objects that are used to describe the
anatomical features and the tumors. Those values are approximate linear
attenuation coefficients per cm at 60 keV (0.416 for bone, 0.210 for
brain, 0.207 for cerebrospinal fluid). The contrast of the small tumors
with their background is 0.003 cm\textsuperscript{-1}. In order to
clearly see the low-contrast details in the interior of the skull,
we use zero (black) to represent the value 0.204 (or anything less)
and 255 (white) to represent 0.21675 or anything more).

For the selected head phantom we generated \textit{parallel projection
data}, in which one \textit{view} comprises estimates of integrals
through the phantom for a set of 693 equally-spaced parallel lines
with a spacing of 0.0376 cm between them. (We chose to simulate parallel
rather than divergent projection data, since the reconstruction by
the method of\cite{SIDK08a} with which we wish to compare the superiorization
approach were performed for us by the authors of\cite{SIDK08a} on
parallel data. Even though contemporary CT scanners use divergent
projection data, results obtained by the use of parallel projection
data are relevant to them, since it is known that the quality of reconstructions
from these two modes of data collection are very similar as long as
the data generations use similar frequencies of sampling of lines
and similar noise characteristics in the estimated integrals for those
lines; see, for example, the reconstructions from divergent and parallel
projection data in figure 5.15 of \cite{HERM09a}.) In calculating
these estimates we take into consideration the effects of photon statistics,
detector width and scatter. Details of how we do this exactly can
be found in Sections 5.5 and 5.9 of \cite{HERM09a}. Briefly, quantum
noise is calculated based on the assumption that approximately 2,000,000
photons enter the head along each ray, detector width is simulated
by using 11 sub-rays along each of which the attenuation is calculated
independently and then combined at the detector, and 5\% of the photons
get counted not by the detector for the ray in question but detectors
for the neighboring rays. For the experiments in this paper, we did
not simulate the poly-energetic nature of the x-ray source. To indicate
what can be achieved in clinical CT, we show in figure \ref{fig:(a)-A-head}(b)
a reconstruction that was made from data comprising of 360 such views
with the reconstruction algorithm known as ART with blob basis functions;
see\cite[Chapter 11]{HERM09a}.

\subsection{Superiorization reconstruction from a few views}

The main reason in the literature for advocating the use of $TV$
as the optimization criterion is that by doing so one can achieve
efficacious reconstructions even from sparsely sampled data. In our
own work\cite{HERM08a} with realistically simulated CT data we found
that this is not always the case and this will be demonstrated again
by the experiments reported in the current paper.

\begin{figure}
\includegraphics[scale=0.5]{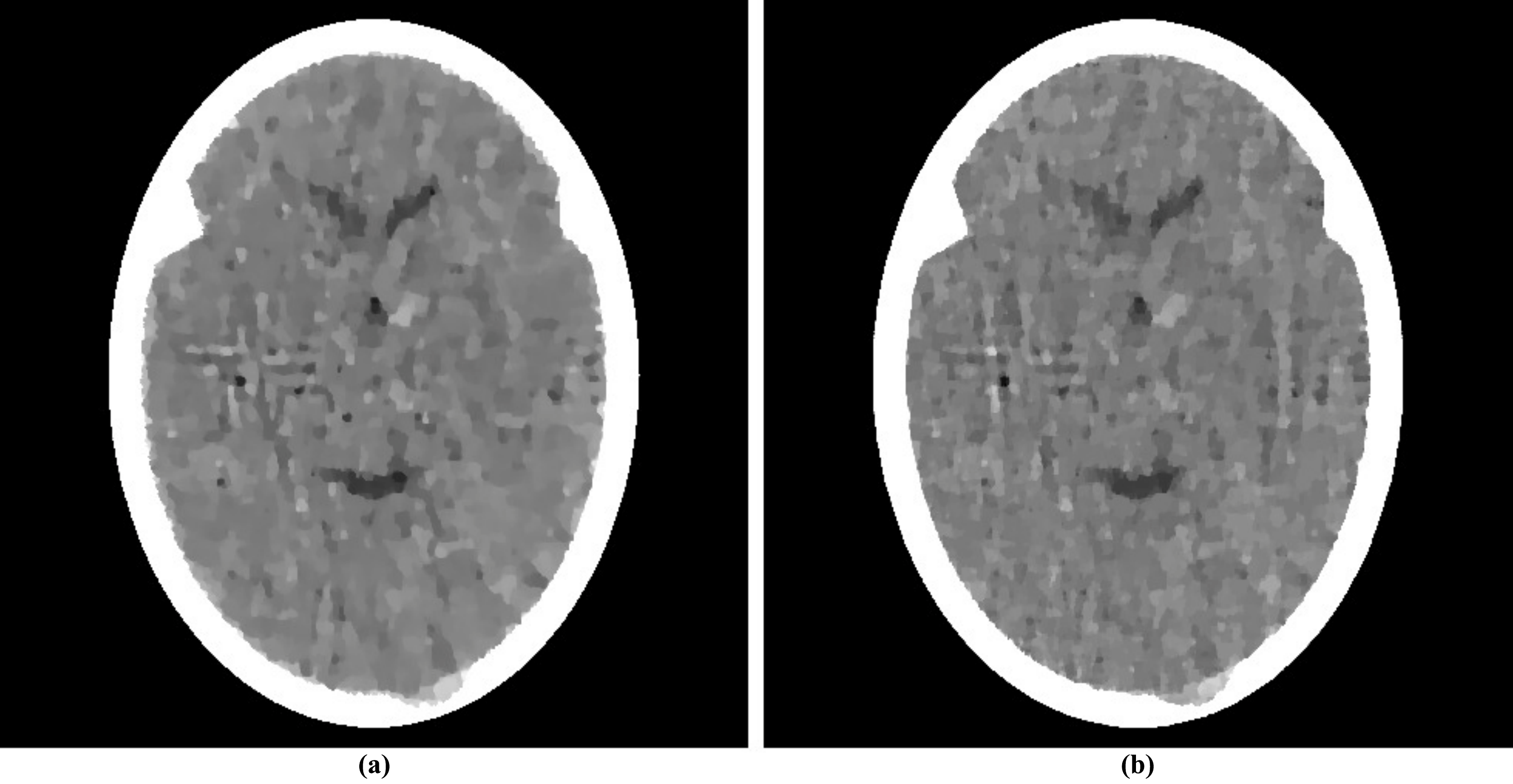}\caption{\label{fig:Reconstructions-using-}Reconstructions using $TV$ as
the optimization criterion from realistically simulated projection
data for 60 views using (a) ASD-POCS and (b) superiorization. As compared
to figure \ref{fig:(a)-A-head}(b), these reconstructions fail in
two ways: they do not show some of the fine details in the phantom
and they present some artifactual variations. The former of these
is a consequence of reconstructing from a much smaller data set than
used for figure \ref{fig:(a)-A-head}(b). The latter is due to using
a very narrow window (13.5 HU) in these displays. Were we to use a
wider display window (e.g., from -429 HU to 429 HU) for the reconstructions
in this figure and in figure \ref{fig:(a)-A-head}(b), the visual
appearance of the resulting images would be nearly indistinguishable.}
\end{figure}

There have appeared in the literature some approaches to $TV$ minimization
that seem to indicate a more efficacious performance for CT than the
one reported in \cite{HERM08a}. One of these is the
Adaptive Steepest Descent Projections Onto Convex Sets (ASD-POCS)
algorithm, which is described in detail in the much-cited paper of
Sidky and Pan\cite{SIDK08a} and whose use has been since reported
in a number of subsequent publications, for example, in \cite{SIDK11a,BIAN10A}.
We note that ASD-POCS was designed with the aim of producing an exact
minimization algorithm, in contrast to our heuristic superiorization
approach. Translating equations (6)-(8) of \cite{SIDK08a} into our
terminology, the aim of ASD-POCS is the following: Given an $\varepsilon\in\mathbb{R}_{+}$,
find an $\varepsilon$-compatible $\boldsymbol{x}\in\Omega=\mathbb{R}_{+}^{J}$
for which $TV(\boldsymbol{x})$ is minimal. (Note that this aim is
a special case of the constrained optimization formulation presented
in (\ref{eq:constrained_minimization}).) In order to test ASD-POCS,
we generated realistic projection data as described in the previous
subsection but for only 60 views at 3 degree increments with the spacing
between the lines for which integrals are estimated set at 0.752 mm.
Thus the number of rays (and hence the number photons put into the
head) in this data set is a twelfth of what it is in the data set
used to produce the reconstruction in figure \ref{fig:(a)-A-head}(b).
A reconstruction from these data was produced for us using ASD-POCS
by the authors of \cite{SIDK08a} (this
ensured that it does not suffer due to our misinterpretation of the
algorithm or from our inappropriate choices of the free parameters),
it is shown in figure \ref{fig:Reconstructions-using-}(a).

Since the image quality of figure \ref{fig:Reconstructions-using-}(a)
is not anywhere near to that of figure \ref{fig:(a)-A-head}(b), we
present here a brief discussion as to why we are showing such images.
Many publications in the recent medical imaging literature have claimed
that medically-efficacious reconstructions can be obtained by the
use of $TV$-minimization from data as sparse as what was used to
produce figure \ref{fig:Reconstructions-using-}(a). (In fact, ASD-POCS
was motivated and used with such an aim in mind\cite{SIDK11a,SIDK08a,BIAN10A}.)
Such publications usually show reconstructions from sparse data as
evidence for the validity of their claims. They can do this because
in their presented illustrations the features that are observable
in the reconstructions are usually much larger and/or of much higher
contrast against their backgrounds than the small ``tumors'' in
figure \ref{fig:(a)-A-head}(a), which are perfectly visible in the
reconstruction in figure \ref{fig:(a)-A-head}(b), but are not detectable
in the reconstruction from sparse data in figure \ref{fig:Reconstructions-using-}(a).
The reason why that reconstruction appears to be unacceptably bad
is that the display window (from 0.204 cm\textsuperscript{-1} linear
attenuation coefficient to 0.21675 cm\textsuperscript{-1} linear
attenuation coefficient) is very narrow; it was selected to enhance
the visibility of the small low-contrast tumors. The width of this
window corresponds to about 13.5 Hounsfield Units (HU). As compared
to this, in their evaluation of sparse-view reconstruction from flat-panel-detector
cone-beam CT, Bian \textit{et al}.\cite{BIAN10A} use what they call
a ``soft-tissue grayscale window'' (also a ``narrow window'')
from -429 HU to 429 HU to display head phantom reconstructions. Using
such a window for our reconstructions shown figures \ref{fig:Reconstructions-using-}(a)
and \ref{fig:(a)-A-head}(b) would result in images that are nearly
indistinguishable from each other. Thus reporting the images using
such a display window is consistent with the claim that a TV-minimizing
reconstruction from a few views is similar in quality to a more traditional
reconstruction from many views. However, our much narrower display
window reveals that this is not really so. We therefore continue using
our much narrower window in what follows, since it clearly reveals
the nature of the reconstructions being compared, warts and all.

While this ASD-POCS reconstruction is not as good as it should be
for diagnostic CT of the brain (due to the sparsity of the data),
it is visually better than the reconstruction using superiorization
from similar data as reported in\cite{HERM08a}. We discuss the reasons
for this in the next subsection. Here we concentrate on examining
whether one can achieve a reconstruction using superiorization that
is as good as that produced by ASD-POCS from the same data.

For this we first need to examine the numerical properties of the
ASD-POCS reconstruction. This reconstruction uses $485\times485$
pixels each of size 0.376 mm by 0.376 mm. This implies that $J=235,225$
and it also determines the components of the vectors $\boldsymbol{a}^{i}\in\mathbb{R}^{J}$
in the precise specification of the problem $S$. The $Res_{S}$,
as defined by (\ref{eq:Res_definition}), of the ASD-POCS reconstruction
is 0.33 and the $TV$, as defined by (\ref{eq:TV}), is 835.

We applied to the same problem $S$ a superiorized version of the
algorithm $\mathbf{R}$ defined by (\ref{eq:R_S_definition}). To
complete the specification of $\mathbf{R}$, we point out that for
the ordering of views we chose the ``efficient'' one that was introduced
in \cite{HERM93a} and is also discussed on page 209 of \cite{HERM09a}.
The choices we made for the superiorization are the following: $\gamma_{\ell}=0.99995^{\ell}$,
$\boldsymbol{\bar{x}}$ is the zero vector and $N=20$. The nonascending
vector was computed by the method described in the paragraph below
(\ref{eq:TV}). Denoting by $R_{S}$ the infinite sequence of points
in $\Omega$ that is produced by the superiorized version of the algorithm
$\mathbf{R}$ when applied to the problem $S$, we chose as our reconstruction
$\boldsymbol{x^{*}}=O\left(S,0.33,R_{S}\right)$. For such a reconstruction
we have, by the definition of $O$, that $Res_{S}\left(\boldsymbol{x^{*}}\right)\leq0.33$;
in other words, the output of the superiorization algorithm is at
least as constraints-compatible with $S$ as the output of ASD-POCS.
From the point of view of $TV$-minimization, our $\boldsymbol{x^{*}}$
is slightly better: $TV\left(\boldsymbol{x^{*}}\right)=$826.

The superiorization reconstruction is displayed in figure \ref{fig:Reconstructions-using-}(b).
Visually it is similar to the reconstruction produced by ASD-POCS.
From the optimization point of view it achieves the desired aim better
than ASD-POCS does, since it results in smaller values for both $Res_{S}$
and for $TV$, even though only slightly.

That the two reconstructions in figure \ref{fig:Reconstructions-using-}
are very similar is not surprising because a comparison of the pseudo-codes
reveals that the ASD-POCS algorithm in
\cite{SIDK08a} is essentially a special case of
the Superiorized Version of Algorithm $\mathbf{P}$, even though
it has been derived from rather different principles. To obtain the
ASD-POCS algorithm from our methodology described here, we would have
to choose an Algebraic Reconstruction Technique (ART; see Chapter
11 of \cite{HERM09a}) as the algorithm that we are superiorizing.
Such a superiorization of ART was reported in the earliest paper on
superiorization \cite{BUTN07a}. For the illustration in our current
paper we decided to superiorize the block-iterative algorithm $\mathbf{R}$
defined by (\ref{eq:R_S_definition}). This illustrates the generality
of the superiorization approach: it is applicable not only to a large
class of constrained optimization problems, but also enables the use
of any of a large class of iterative algorithms designed to produce
a constraints-compatible solutions. A recent publication aimed at
producing an exact $TV$-minimizing algorithm based on the block-iterative
approach is \cite{DEFR11a}.

\subsection{\label{sub:Effects-of-variations}Effects of variations in the reconstruction
approach}

The reconstruction in figure \ref{fig:Reconstructions-using-}(a)
produced by ASD-POCS definitely ``looks better'' than a reconstruction
in \cite{HERM08a}, which was obtained using
superiorization from similar data. Since, as discussed in the last
paragraph of the previous subsection, the ASD-POCS
algorithm in \cite{SIDK08a} can be obtained as a special case of
superiorization, it must be that some of the choices made in the
details of the implementations are responsible for the visual differences.
An analysis of the implementational details adopted by the two approaches
revealed several differences. After removing these differences, the
superiorization approach produced the image in figure \ref{fig:Reconstructions-using-}(b),
which is very similar to the reconstruction produced by ASD-POCS.
We now list the implementational choices that were made for superiorization
to make its performance match that of the reported implementation
of ASD-POCS.

\begin{figure}
\includegraphics[scale=0.5]{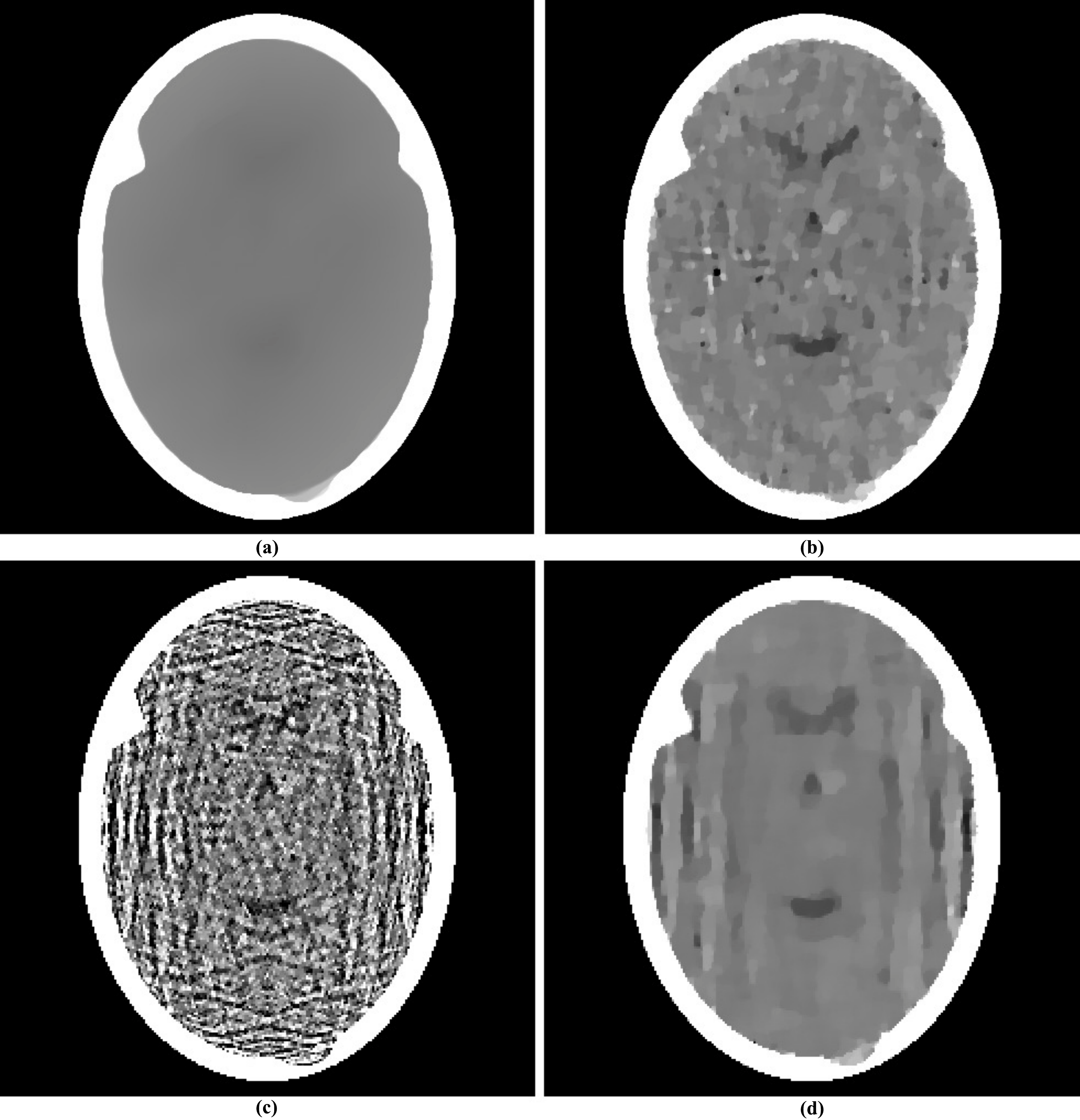}
\caption{
\label{fig:Reconstructions-produced-by}Reconstructions produced by
varying some of the parameters in the algorithm that produced figure
\ref{fig:Reconstructions-using-}(b). (a) Changing the termination
criterion form $\varepsilon=0.33$ to $\varepsilon=0.91.$ (b) Changing
the value of $N$ from 20 to 1. (c) Reconstructing with pixel size
0.752 mm by 0.752 mm instead of 0.376 mm by 0.376 mm. (d) Reconstructing
with all the three changes of (a)-(c).
}
\end{figure}

One implementational difference is in the stopping-rule of the iterative
algorithm; that is, the choice of $\varepsilon$ in determining the
output $O\left(S,\varepsilon,R_{S}\right)$. Since the data are noisy,
the phantom itself does not match the data exactly. In previously
reported implementations of superiorization it was assumed that the
iterative process should terminate when an image is obtained that
is approximately as constraints-compatible as the phantom; in the
case of the phantom and the projections data on which we report here
the value of $Res_{S}$ for the phantom
is approximately 0.91, which is larger than its value (0.33) for the
reconstruction produced by ASD-POCS. The output $O\left(S,0.91,R_{S}\right)$
is shown in figure \ref{fig:Reconstructions-produced-by}(a). This
is a wonderfully smooth reconstruction, its $TV$ value is only 771.
However this smoothness comes at a price: we loose not only the ability
to detect the large tumor, but we cannot even see anatomic features
(such as the ventricular cavities) inside the brain. So it appears
that, in order to see medically-relevant features in the brain, \textit{over-fitting}
(in the sense of producing a reconstruction from noisy data that is
more constraints-compatible than the phantom) is desirable.

In the implementations that produced previously reported reconstructions
by superiorization, the number $N$ in the Superiorized
Version of Algorithm
$\mathbf{P}$ was always chosen to be 1. It is possible that this
is the wrong choice, making only this change to what lead to the reconstruction
in figure \ref{fig:Reconstructions-using-}(b)
results in the reconstruction shown in figure \ref{fig:Reconstructions-produced-by}(b).
That image appears similar to the image in figure
\ref{fig:Reconstructions-using-}(b), but it has a higher $TV$ value,
namely 832, which is still very slightly lower than that of the ASD-POCS
reconstruction. The choice $N=20$ was based on the desire to maintain
consistency with what has been practiced using ASD-POCS, see page
4790 of \cite{SIDK08a}. It appears that
in the context of our paper the additional computing cost due to choosing
$N$ to be 20 rather than 1 is not really justified. (We note that
if $\boldsymbol{d}$ is selected using subgradients as discussed in
the paragraph following (\ref{eq:nonascending}) and thus $\boldsymbol{d}$
is not guaranteed to be a nonascending vector for the $TV$ function,
then the choice of 20 rather than 1 for $N$ results in a considerable
improvement. However, an even greater improvement is achieved even
with $N=1$ by selecting $\boldsymbol{d}$ as recommended in this
paper.)

Another important difference between the
ASD-POCS implementation and the previous implementations of the superiorization
approach is the size of the pixels in the reconstructions. For the
ASD-POCS reconstruction this was selected to be 0.376 mm by 0.376
mm. In previously reported reconstructions by superiorization it was
assumed that the edge of a pixel should be the same as the distance
between the parallel lines along which the data are collected; that
is, 0.752 mm for our problem $S$. This assumption proved to be false.
$TV$-minimization takes care of undesirable artifacts that may otherwise
arise due to the smaller pixels and this leads to a visual improvement.
A superiorizing reconstruction with the larger pixels, using $\varepsilon=0.33$
and $N=20$, is shown in figure \ref{fig:Reconstructions-produced-by}(c).
(We note that the use of smaller pixels during iterative x-ray CT
reconstructions was also suggested in\cite{ZBIJ04a}. However, that
approach is quite different from what is presented here: its final
result uses larger pixels whose values are obtained by averaging assemblies
of values provided by the iterative process to the smaller pixels.
There is no such downsampling in our approach, our final result is
presented using the smaller pixels. Its smoothness is due to reduction
of TV by the superiorization approach rather than to averaging pixel
values in a denser digitization.)

Combining the use of the larger pixels with $\varepsilon=0.91$ and
$N=1$ results in the reconstruction shown in figure \ref{fig:Reconstructions-produced-by}(d).
This reconstruction, for which the superiorization options were selected
according to what was done in\cite{HERM08a}, is visually inferior
to those shown in our figure \ref{fig:Reconstructions-using-}. The
reconstructions displayed in figure \ref{fig:Reconstructions-produced-by}
also illustrate another important point, namely that even though the
mathematical results discussed in this paper are valid for a large
range of choices of the parameters in the superiorization algorithms,
for medical efficacy of the reconstructions attention has to be paid
to these choices since they can have a drastic effect on the quality
of the reconstruction.

\begin{figure}
\includegraphics[scale=0.5]{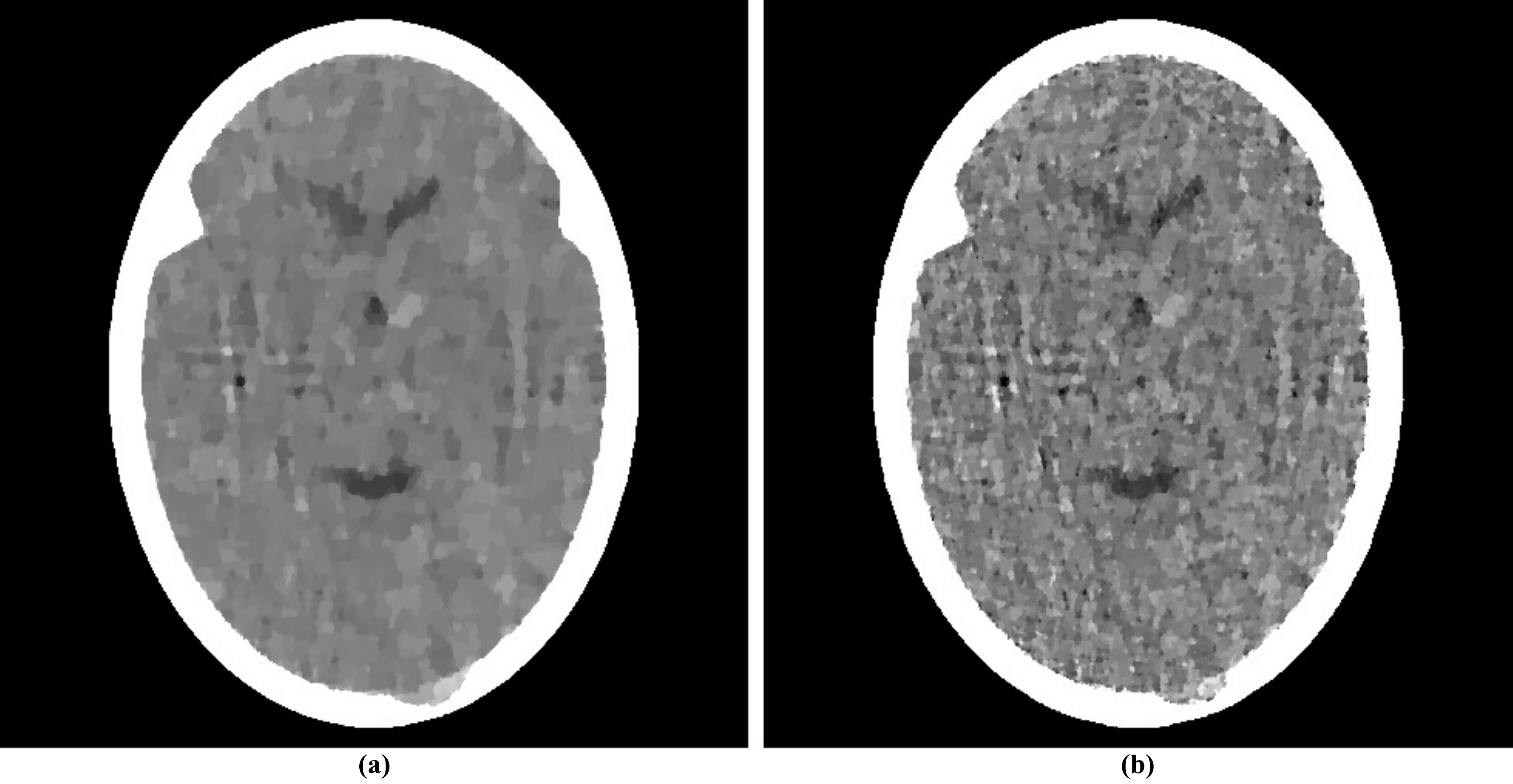}\caption{
\label{fig:Reconatructions-by-variants}Reconstructions by variations
that do not fit into the framework within which the previously shown
reconstructions were produced. (a) Not using nonnegativity in the
algorithm. (b) Interleaving perturbations with blocks.
}
\end{figure}

It has been mentioned in Subsection \ref{sub:Algorithms-and-outputs}
that except for the presence of $\mathbf{Q}$
in (\ref{eq:R_S_definition}), which enforces nonnegativity of the
components, $\mathbf{R}$ is identical to the algorithm used and illustrated
in \cite{HERM08a}. It is known that CT reconstruction of the brain
from many views does not suffer from ignoring the fact that the components
of the $\boldsymbol{x}$, which represent linear attenuation coefficients,
should be nonnegative; as is illustrated in figure \ref{fig:(a)-A-head}(b).
This remains so when reconstructing from a few views using the method
and data that we have been discussing: if we do everything in exactly
the same way as was done to obtain the reconstruction with $TV$ value
826 that is shown in our figure \ref{fig:Reconstructions-using-}(b)
but remove $\mathbf{Q}$ from (\ref{eq:R_S_definition}), then we
obtain a reconstruction in figure \ref{fig:Reconatructions-by-variants}(a)
whose $TV$ value is 829.

Another variation that deserves discussion, because it has been suggested
in the literature \cite{PENF10a}, is one
that does not come about by making choices for the general approach
of the Superiorized Version of Algorithm
$\mathbf{P}$ but rather by changing the nature of the approach. The
variation in question is not applicable in general, but can be applied
to the special case when the algorithm to be superiorized is the $\mathbf{R}$
defined by (\ref{eq:R_S_definition}). It was suggested as an improvement
to the approach presented above with the choice $N=1$. The idea was
based on recognizing the block-iterative nature of the algorithmic
operator $\mathbf{R}_{S}$ in (\ref{eq:R_S_definition}) and intermingling
the perturbation steps of lines (vii)-(xvii) of the Superiorized Version
of Algorithm $\mathbf{R}$ with the projection steps $\mathbf{B}_{S_{1}},\ldots,\mathbf{B}_{S_{W}}$
of (\ref{eq:R_S_definition}). It was reported in \cite{PENF10a}
that doing this is advantageous to using the Superiorized Version
of Algorithm $\mathbf{R}$. However, when we applied the variation
of the Superiorized Version of Algorithm $\mathbf{R}$ that is proposed
in \cite{PENF10a} to the problem $S$ that we have been using in
this section, we ended up with the reconstruction in figure \ref{fig:Reconatructions-by-variants}(b)
whose $TV$ value is 920. This is not as good as what was obtained
using the version of the algorithm that produced the reconstruction
in figure \ref{fig:Reconstructions-using-}(b). We conclude that
the variation suggested by \cite{PENF10a},
which does not fit into the theory of our paper, does not have an
advantage over what we are proposing here, at least for the problem
$S$ that we have been discussing in this section. We conjecture that
the improvement reported in \cite{PENF10a} is due to selecting $\boldsymbol{d}$
using subgradients as discussed in the paragraph following (\ref{eq:nonascending})
and, as discussed earlier, such an improvement is not obtained if
$\boldsymbol{d}$ is selected by the more appropriate method recommended
in this paper.

\section{Discussion and Conclusions}

Constrained optimization is an often-used tool in medical physics.
The methodology of superiorization is a heuristic (as opposed to exact)
approach to constrained optimization.

Although the idea of superiorization was introduced in 2007 and its
practical use has been demonstrated in several publications since,
this paper is the first to provide a solid mathematical foundation
to superiorization as applied to the noisy problems of the real world.
These foundations include a precise definition of constraints-compatibility,
the concept of a strongly perturbation resilient algorithm, simple
conditions that ensure that an algorithm is strongly perturbation
resilient, the superiorized version of an algorithm and the showing
that the superiorized version of a strongly perturbation resilient
algorithm produces outputs that are essentially as constraints-compatible
as those produced by the original version but are likely to have a
smaller value of the chosen optimization criterion.

The approach is very general. For any iterative algorithm $\mathbf{P}$
and for any optimization criterion $\phi$ for which we know how to
produce nonascending vectors, the pseudocode given in Subsection \ref{sub:Superiorized-version-of}
automatically provides the version of $\mathbf{P}$ that is superiorized
for $\phi$.

We demonstrated superiorization for tomography when total variation
is used as the optimization criterion. In particular, we illustrated
on a particular tomography problem that, in spite of its generality,
superiorization produced a reconstruction that is as good as (from
the points of view of constraints-compatibility and $TV$-minimization)
what was obtained by the ASD-POCS algorithm that was specially designed
for $TV$-minimization in tomography.

\section*{Acknowledgments}

The detailed and penetrating comments of three reviewers and the editors
helped us to improve this paper in a significant way. We thank Prof.
Xiaochuan Pan and his coworkers from the University of Chicago for
providing us with the reconstruction from our data using their implementation
of their ASD-POCS algorithm. Our work is supported by the National
Science Foundation award number DMS-1114901, the United States-Israel
Binational Science Foundation (BSF) grant number 200912, and the US
Department of Army award number W81XWH-10-1-0170.

\section*{Appendix}

\subsection*{Conditions for strong perturbation resilience}

\paragraph*{Theorem 1.}

Let $\mathbf{P}$ be an algorithm for a problem structure $\left\langle \mathbb{T},\mathcal{P}r\right\rangle $
such that, for all $T\in\mathbb{T}$, $\mathbf{P}$ is boundedly convergent
for $T$, $\mathcal{P}r_{T}:\Omega\rightarrow\mathbb{R}$ is uniformly
continuous and $\mathbf{P}_{T}:\Delta\rightarrow\Omega$ is nonexpansive.
Then $\mathbf{P}$ is strongly perturbation resilient.

\paragraph*{Proof.}

We first show that there exists an $\varepsilon\in\mathbb{R}_{+}$
such that $O\left(T,\varepsilon,\left(\left(\mathbf{P}_{T}\right)^{k}\boldsymbol{x}\right)_{k=0}^{\infty}\right)$
is defined for every $\boldsymbol{x}\in\Omega$. Under the assumptions
of the theorem, let $\gamma\in\mathbb{R}_{+}$ be such that $\mathcal{P}r_{T}\left(\boldsymbol{y}\left(\boldsymbol{x}\right)\right)\leq\gamma$,
for every $\boldsymbol{x}\in\Omega$. We prove that $O\left(T,2\gamma,\left(\left(\mathbf{P}_{T}\right)^{k}\boldsymbol{x}\right)_{k=0}^{\infty}\right)$
is defined for every $\boldsymbol{x}\in\Omega$ as follows. Select
a particular $\boldsymbol{x}\in\Omega$. By uniform continuity of
$\mathcal{P}r_{T}$, there exists a $\delta>0$, such that $\left|\mathcal{P}r_{T}\left(\boldsymbol{z}\right)-\mathcal{P}r_{T}\left(\boldsymbol{y}\left(\boldsymbol{x}\right)\right)\right|\leq\gamma$,
for any $\boldsymbol{z}\in\Omega$ for which $\left\Vert \boldsymbol{z}-\boldsymbol{y}\left(\boldsymbol{x}\right)\right\Vert \leq\delta$.
Since $\mathbf{P}$ is convergent for $T$, there exists a nonnegative
integer $K$, such that $\left\Vert \left(\mathbf{P}_{T}\right)^{K}\boldsymbol{x}-\boldsymbol{y}\left(\boldsymbol{x}\right)\right\Vert \leq\delta$.
It follows that

\begin{equation}
\begin{array}{rcl}
\left|\mathcal{P}r_{T}\left(\left(\mathbf{P}_{T}\right)^{K}\boldsymbol{x}\right)\right| & \leq & \left|\mathcal{P}r_{T}\left(\left(\mathbf{P}_{T}\right)^{K}\boldsymbol{x}\right)-\mathcal{P}r_{T}\left(\boldsymbol{y}\left(\boldsymbol{x}\right)\right)\right|+\left|\mathcal{P}r_{T}\left(\boldsymbol{y}\left(\boldsymbol{x}\right)\right)\right|\\
 & \leq & 2\gamma.
\end{array}\label{eq:convergence}
\end{equation}

Now let $T\in\mathbb{T}$ and $\varepsilon\in\mathbb{R}_{+}$ be such
that $O\left(T,\varepsilon,\left(\left(\mathbf{P}_{T}\right)^{k}\boldsymbol{x}\right)_{k=0}^{\infty}\right)$
is defined for every $\boldsymbol{x}\in\Omega$. To prove the theorem,
we need to show that $O\left(T,\varepsilon',R\right)$ is defined
for every $\varepsilon'>\varepsilon$ and for every sequence $R=\left(\boldsymbol{x}^{k}\right)_{k=0}^{\infty}$
of points in $\Omega$ for which, for all $k\geq0$, (\ref{eq:perturbations})
is satisfied for bounded perturbations $\beta_{k}\boldsymbol{v}^{k}$
. Let $\varepsilon'$ and $R$ satisfy the conditions of the previous
sentence.

For $k\geq0$, we have, due to the nonexpansiveness of $\mathbf{P}_{T}$,
that 
\begin{equation}
\left\Vert \boldsymbol{x}^{k+1}-\mathbf{P}_{T}\boldsymbol{x}^{k}\right\Vert =\left\Vert \mathbf{P}_{T}\left(\boldsymbol{x}^{k}+\beta_{k}\boldsymbol{v}^{k}\right)-\mathbf{P}_{T}\boldsymbol{x}^{k}\right\Vert \leq\left\Vert \beta_{k}\boldsymbol{v}^{k}\right\Vert .\label{eq:nonexpansiveness}
\end{equation}
Denote $\left\Vert \beta_{k}\boldsymbol{v}^{k}\right\Vert $ by $r_{k}$.
Clearly, $r_{k}\in\mathbb{R}_{+}$ and it follows from the definition
of bounded perturbations that ${\displaystyle \sum\limits _{k=0}^{\infty}}r_{k}\,<\infty$.

We next prove by induction that, for every pair of nonnegative integers
$k$ and $i$, 
\begin{linenomath}
\begin{equation}
\left\Vert \boldsymbol{x}^{k+i}-\left(\mathbf{P}_{T}\right)^{i}\boldsymbol{x}^{k}\right\Vert \leq{\displaystyle \sum_{j=k}^{k+i-1}r_{j}}.\label{eq:Lemma}
\end{equation}
\end{linenomath}
Let $k$ be an arbitrary nonnegative integer. If $i=0$, then the
value is zero on both sides of the inequality and hence (\ref{eq:Lemma})
holds. Now assume that (\ref{eq:Lemma}) holds for an integer $i\geq0$.
Then, by (\ref{eq:nonexpansiveness}) and the nonexpansiveness of
$\mathbf{P}_{T}$, 
\begin{linenomath}
\begin{equation}
\begin{array}{rcl}
\left\Vert \boldsymbol{x}^{k+i+1}-\left(\mathbf{P}_{T}\right)^{i+1}\boldsymbol{x}^{k}\right\Vert  & \leq & \left\Vert \boldsymbol{x}^{k+i+1}-\mathbf{P}_{T}\boldsymbol{x}^{k+i}\right\Vert \\
 &  & +\left\Vert \mathbf{P}_{T}\boldsymbol{x}^{k+i}-\left(\mathbf{P}_{T}\right)^{i+1}\boldsymbol{x}^{k}\right\Vert \\
 & \leq & r_{k+i}+\left\Vert \boldsymbol{x}^{k+i}-\left(\mathbf{P}_{T}\right)^{i}\boldsymbol{x}^{k}\right\Vert \\
 & \leq & r_{k+i}+{\displaystyle \sum_{j=k}^{k+i-1}r_{j}}\\
 & = & {\displaystyle \sum_{j=k}^{k+i}r_{j}},
\end{array}\label{eq:induction}
\end{equation}
\end{linenomath}
which completes our inductive proof. A consequence of (\ref{eq:Lemma})
is that, for every pair of nonnegative integers $k$ and $i$, 
\begin{linenomath}
\begin{equation}
\left\Vert \boldsymbol{x}^{k+i}-\left(\mathbf{P}_{T}\right)^{i}\boldsymbol{x}^{k}\right\Vert \leq{\displaystyle \sum_{j=k}^{\infty}r_{j}}.\label{eq:Lemma-1}
\end{equation}
\end{linenomath}
Due to the summability of the nonnegative sequence $\left(r_{k}\right)_{k=0}^{\infty}$,
the right-hand side (and hence the left-hand side) of this inequality
gets arbitrarily close to zero as $k$ increases.

Since $\mathcal{P}r_{T}$ is uniformly continuous, there exists a
$\delta$ such that, for all $\boldsymbol{x},\boldsymbol{y}\in\Omega$,
$\left|\mathcal{P}r_{T}(\boldsymbol{x})-\mathcal{P}r_{T}(\boldsymbol{y})\right|\leq\varepsilon'-\varepsilon$
provided that $\left\Vert \boldsymbol{x}-\boldsymbol{y}\right\Vert \leq\delta$.
Select a $k$ so that $\sum_{j=k}^{\infty}r_{j}\leq\delta$. By the
assumption that $O\left(T,\varepsilon,\left(\left(\mathbf{P}_{T}\right)^{k}\boldsymbol{x}\right)_{k=0}^{\infty}\right)$
is defined for every $\boldsymbol{x}\in\Omega$, there exists a nonnegative
integer $i$ for which $\mathcal{P}r\left(\left(\mathbf{P}_{T}\right)^{i}\boldsymbol{x}^{k}\right)\leq\varepsilon$.
From (\ref{eq:Lemma-1}) we have, for this $k$ and $i$, that $\left\Vert \boldsymbol{x}^{k+i}-\left(\mathbf{P}_{T}\right)^{i}\boldsymbol{x}^{k}\right\Vert \leq\delta$
and, hence, 
\begin{linenomath}
\begin{equation}
\begin{array}{rcl}
\left|\mathcal{P}r_{T}(\boldsymbol{x}^{k+i})\right| & \leq & \left|\mathcal{P}r_{T}(\boldsymbol{x}^{k+i})-\mathcal{P}r_{T}\left(\left(\mathbf{P}_{T}\right)^{i}\boldsymbol{x}^{k}\right)\right|\\
 &  & +\left|\mathcal{P}r_{T}\left(\left(\mathbf{P}_{T}\right)^{i}\boldsymbol{x}^{k}\right)\right|\\
 & \leq & \left(\varepsilon'-\varepsilon\right)+\varepsilon=\varepsilon',
\end{array}\label{eq:consistency}
\end{equation}
\end{linenomath}
proving that $O\left(T,\varepsilon',R\right)$ is defined. $\square$

\subsection*{Nonascending vectors for convex functions}

\paragraph*{Theorem 2.}

Let $\phi:\mathbb{R}^{J}\rightarrow\mathbb{R}$ be a convex function
and let $\boldsymbol{x}\in\mathbb{R}^{J}$. Let $\boldsymbol{g}\in\mathbb{R}^{J}$
satisfy the property: For 1$\leq j\leq J$, if the $j$th component
$g_{j}$ of $\boldsymbol{g}$ is not zero, then the partial derivative
$\frac{\partial\phi}{\partial x_{j}}(\boldsymbol{x})$ of $\phi$
at $\boldsymbol{x}$ exists and its value is $g_{j}$. Define $\boldsymbol{d}$
to be the zero vector if $\left\Vert \boldsymbol{g}\right\Vert =0$
and to be $-\boldsymbol{g}/\left\Vert \boldsymbol{g}\right\Vert $
otherwise. Then $\boldsymbol{d}$ is a nonascending vector for $\phi$
at $\boldsymbol{x}$.

\paragraph{Proof. }

The theorem is trivially true if $\left\Vert \boldsymbol{g}\right\Vert =0$,
so we assume that this is not the case. We denote by $I$ the nonempty
set of those indices $j$ for which $g_{j}\neq0$.

For $1\leq j\leq J$, let $s_{j}$ be $\nicefrac{g_{j}}{\left|g_{j}\right|}$
for $j\in I$ and be 0 otherwise, and let $\boldsymbol{e}^{j}\in\mathbb{R}^{J}$
be the vector all of whose components are zero except for the $j$th,
which is one. Then, for $1\leq j\leq J$, there exists a $\delta_{j}>0$
such that, for $0\leq\lambda_{j}\leq\delta_{j}$, 
\begin{linenomath}
\begin{equation}
\phi\left(\boldsymbol{x}-\lambda_{j}s_{j}\boldsymbol{e}^{j}\right)\leq\phi\left(\boldsymbol{x}\right).\label{eq:descending}
\end{equation}
\end{linenomath}
This is obvious if $s_{j}=0$. Otherwise, $\frac{\partial\phi}{\partial x_{j}}(\boldsymbol{x})$
exists and indicates $\phi$ increases at $\boldsymbol{x}$ if $s_{j}=1$
or that $\phi$ decreases at $\boldsymbol{x}$ if $s_{j}=-1$. The
existence of the desired $\delta_{j}$ can be derived from the standard
definition of the partial derivative as a limit.

We define $\delta>0$ by 
\begin{linenomath}
\begin{equation}
\delta=\frac{\left\Vert \boldsymbol{g}\right\Vert }{J}\min_{j\in I}\left\{ \frac{\delta_{j}}{\left|g_{j}\right|}\right\} .\label{eq:delata_definition}
\end{equation}
\end{linenomath}
Then we have that, for $0\leq\lambda\leq\delta$, 
\begin{linenomath}
\begin{equation}
\begin{array}{rcl}
\phi\left(\boldsymbol{x}+\lambda\boldsymbol{d}\right) & = & \phi\left(\boldsymbol{x}-\lambda{\displaystyle \sum_{j=1}^{J}\frac{\left|g_{j}\right|}{\left\Vert \boldsymbol{g}\right\Vert }s_{j}}\boldsymbol{e}^{j}\right)\\
 & = & \phi\left({\displaystyle \sum_{j=1}^{J}\frac{1}{J}\left(\boldsymbol{x}-\lambda J\frac{\left|g_{j}\right|}{\left\Vert \boldsymbol{g}\right\Vert }s_{j}\boldsymbol{e}^{j}\right)}\right)\\
 & \leq & {\displaystyle \frac{1}{J}{\displaystyle \sum_{j=1}^{J}\phi\left(\boldsymbol{x}-\lambda J\frac{\left|g_{j}\right|}{\left\Vert \boldsymbol{g}\right\Vert }s_{j}\boldsymbol{e}^{j}\right)}}\\
 & \leq & {\displaystyle \frac{1}{J}\sum_{j=1}^{J}}\phi\left(\boldsymbol{x}\right)\\
 & = & \phi\left(\boldsymbol{x}\right).
\end{array}\label{eq:descent_proof}
\end{equation}
\end{linenomath}
The first inequality above follows from the convexity of $\phi$ and
the second one follows from (\ref{eq:descending}), with $\lambda_{j}$
defined to be $\lambda J\frac{\left|g_{j}\right|}{\left\Vert \boldsymbol{g}\right\Vert }$,
combined with (\ref{eq:delata_definition}). Thus $\boldsymbol{d}$
is a nonascending vector for $\phi$ at $\boldsymbol{x}$. $\square$\bigskip{}

\end{document}